\documentclass[a4paper]{amsart}
\usepackage{graphicx} 
\usepackage{amsmath}
\usepackage{amsfonts}
\usepackage{mathrsfs}
\usepackage{graphicx}
\usepackage{amsmath,amsthm,amssymb,amscd}
\usepackage{color}
\usepackage{xcolor}
\usepackage{enumerate}
\usepackage[T1]{fontenc}
\usepackage{tikz-cd}
\usetikzlibrary{calc}
\newtheorem{thm}{Theorem}[section]
\newtheorem{lemma}[thm]{Lemma}
\newtheorem{proposition}[thm]{Proposition}
\newtheorem{defn}[thm] {Definition}
\newtheorem{exmp}[thm]{Example}
\newtheorem{cor}[thm]{Corollary}
\newtheorem{rem} [thm]{Remark}

\newtheorem{que}[thm]{Question}

\DeclareMathOperator{\cl}{cl}

\title{Universal Realization of Milnor Attractors in Low Dimensions}
\date{}

\author[J. Hant\'akov\'a]{Jana Hant\'akov\'a$^{*}$}
\author[M. Ml\'ichov\'a]{Michaela Ml\'ichov\'a}

\thanks{$^{*}$Corresponding author}

\address[J. Hant\'akov\'a]{Silesian University in Opava, Mathematical Institute,
Na Rybn\'i\v cku 626/1, 746 01 Opava, Czech Republic}
\email{jana.hantakova@math.slu.cz}

\address[M. Ml\'ichov\'a]{Silesian University in Opava, Mathematical Institute,
Na Rybn\'i\v cku 626/1, 746 01 Opava, Czech Republic}
\email{michaela.mlichova@math.slu.cz}

\begin{document}

\begin{abstract}
We study the realization problem for Milnor attractors on compact metric spaces and establish a general sufficient criterion for universal realization. In particular, we prove that every smooth compact one-dimensional manifold, with or without boundary, and every smooth closed two-dimensional manifold has the property that each nonempty closed subset is the Milnor global attractor of a continuous self-map. If the manifold is connected, the realizing map may be chosen to be a homeomorphism. We also investigate the realization problem on selected
one-dimensional continua and exhibit spaces for which universal
realization fails.


\end{abstract}

\maketitle

\section{Introduction}

A discrete dynamical system is an ordered pair $(X,f)$, where $X$ is a compact
metric space and $f\colon X\to X$ is continuous. A central problem in
topological dynamics is to understand the asymptotic behavior of points under
iteration of $f$. For a point $x\in X$, this behavior is described by its
$\omega$-limit set, the collection of all points that the orbit of $x$
approaches arbitrarily closely as time goes to infinity.

One of the most important concepts in this context is that of an
\emph{attractor}, intended to represent the asymptotic behavior of a large set
of initial conditions. The meaning of ``large'' has been formalized in
different ways, leading to several nonequivalent notions of attractor. A
measure-theoretic approach was introduced by Milnor~\cite{Milnor2}, who defined
a \emph{Milnor attractor} as the smallest closed set containing the
$\omega$-limit sets of almost all points, with respect to a given probability
measure. This definition allows for exceptional sets of measure zero and is
particularly well suited to applications in which initial conditions cannot be
specified with infinite precision.

Despite its conceptual simplicity, Milnor's definition raises fundamental
questions about the flexibility and rigidity of attractors. For Conley
attractors and several other classical notions of attractor, such flexibility
is known to be limited \cite{Conley1978}. By contrast,
even in simple settings, it is not well understood how strongly the structure
of a Milnor attractor is constrained by the underlying space and the dynamics.
This motivates the following basic problem.

\begin{que}
For which compact metric spaces $X$ does every nonempty closed subset
$A\subseteq X$ occur as a Milnor attractor of some continuous map $f:X\to X$?
\end{que}

Whenever this property holds, we say that $X$ is an \emph{MAF--space} (Milnor
Attractors = family of closed sets). This terminology formalizes the
universality property that every nonempty closed subset can arise as an
observable attractor. It is easy to verify that the interval $[0,1]$ is an
MAF--space. The main purpose of this paper is to determine how far this
flexibility extends in low dimensions.

Our first main result shows that smooth manifolds of low dimension exhibit
maximal freedom.

\begin{thm}\label{thm:intro}
Every smooth compact one-dimensional manifold, with or without boundary,
and every smooth closed two-dimensional manifold is an MAF--space.
Moreover, if \(X\) is connected, then for every nonempty closed subset
\(A\subset X\), the realizing map \(f\) can be chosen to be a homeomorphism.
\end{thm}

\medskip
The one-dimensional case of Theorem~\ref{thm:intro} is
Theorem~\ref{thm:1D}, and the two-dimensional case is
Theorem~\ref{thm:2D}. In both cases the proof is constructive. Given a
closed set $A$, we explicitly construct maps that contract each connected
component of the complement $X\setminus A$ toward its boundary while fixing
$A$ pointwise. This mechanism is formalized via the notion of a
\emph{root--contracting decomposition}, which provides a local method for
forcing almost every orbit to converge to $A$.

As consequences of the surface construction, we obtain two further
results. First, using the Oxtoby--Ulam theorem, the two-dimensional
realization result extends from smooth surfaces with area measure to
closed topological surfaces equipped with any nonatomic full-support
Borel probability measure. Second, the realizing homeomorphism may be
chosen so that the prescribed set \(A\) is simultaneously the Milnor
global attractor and the generic limit set.

In Section~\ref{sec:beyond} we discuss further examples beyond smooth
manifolds and collect a number of open questions that this abstract
mechanism suggests, including the extent to which the root--contracting
approach applies to more general continua.

These flexibility results stand in sharp contrast with what is known about
generic dynamics. For interval maps, recent work of Foryś--Krawiec,
Hantáková, Kowalewski, and Oprocha~\cite{FHKO} shows that for a residual
subset of $C([0,1])$, the Milnor global attractor of a typical map is a
Cantor set of zero Hausdorff dimension, whose topological structure is
highly rigid and cannot be prescribed. Our results show that this generic rigidity is compatible with considerable freedom of realization: every nonempty closed subset of the manifolds
considered here can occur as a Milnor global attractor. The realizing maps
fix the prescribed set pointwise. Thus the topology of the prescribed set
alone imposes no obstruction to realization in these settings; prescribing
nontrivial dynamics on that set remains a separate problem.

\subsection*{Organization of the paper}
The paper is organized as follows. Section~2 introduces the basic
definitions used throughout the paper. In Section~3 we introduce the
notion of root--contracting decompositions and establish general criteria
ensuring the existence of Milnor attractors. Section~4 applies these tools
to smooth compact manifolds of dimension one, and Section~5 to smooth
compact manifolds of dimension two. Section~6 discusses examples beyond
manifolds and collects open questions.

\section{Preliminaries}

By a manifold we mean a smooth, second-countable Hausdorff manifold,
possibly disconnected. One-dimensional manifolds are allowed to have
boundary, whereas two-dimensional manifolds are assumed to have empty
boundary unless stated otherwise.

Let \((X,d)\) be a compact metric space and \(f:X\to X\) a continuous
map. For \(x\in X\), the \emph{orbit} is
\[
\operatorname{Orb}(x)=\{f^n(x):n\geq0\},
\]
and the \emph{\(\omega\)-limit set} is
\[
\omega(x)=\{y\in X:\exists\,n_k\to\infty
\text{ with }f^{n_k}(x)\to y\},
\]
which is always closed and \(f\)-invariant.

Throughout, every compact metric space is understood to be equipped with
a fixed Borel probability measure $\lambda$, called the reference measure.

\subsection{Attractors}

For any nonempty set \(U\subset X\), the basin of attraction of \(U\) is
\[
B(U)=\{x\in X:\omega(x)\subseteq U\}.
\]

\begin{defn}[Milnor global attractor]\label{def:ga}
The Milnor global attractor (or likely limit set) of \(f\) is the
smallest closed set containing the \(\omega\)-limit sets of
\(\lambda\)-almost all orbits.
\end{defn}
Existence and uniqueness of this smallest closed set are classical; see
\cite{Milnor2}.
\begin{defn}[MAF-space]
A compact space \(X\) is an \emph{MAF-space} if every nonempty closed
subset \(A\subseteq X\) is the Milnor global attractor of some
continuous self-map \(f:X\to X\).
\end{defn}
The MAF property is understood relative to the fixed reference measure
\(\lambda\), which we suppress from the notation.
\begin{defn}[Generic limit set]\label{def:generic-limit-set}
Let $X$ be a compact metric space and $f:X\to X$ continuous. The
\emph{generic limit set} $\Gamma(f)$ is the smallest closed set $A\subset
X$ such that $\{x\in X:\omega(x)\subset A\}$ is residual in $X$.
\end{defn}

There are two natural meanings of "a large set of initial conditions": full
measure and residual; the Milnor global attractor $\Lambda(f)$ uses the
former and the generic limit set $\Gamma(f)$ the latter, and in general
$\Lambda(f)\neq\Gamma(f)$ by \cite{Milnor2}.
\subsection{Choice of reference measure}

Milnor attractors are defined in terms of the measure of their basins, so a
reference measure on the ambient space must be fixed. Throughout the paper
we use natural geometric measures associated with the spaces under
consideration; these measures are not assumed to be invariant under the
dynamics. Since multiplying a finite measure by a positive constant does not change
its null sets or full-measure sets, we normalize all finite reference
measures to have total mass one.

For smooth compact manifolds we use the Lebesgue measure induced by an
arbitrary Riemannian metric. Since any two such measures are mutually
absolutely continuous, the notions of Milnor attractor and Milnor global
attractor are independent of the choice of Riemannian metric.

In the examples of Section~\ref{sec:beyond}, we use finite Borel measures
equivalent to the appropriate Hausdorff measure. On finite graphs this is
simply \(H^1\), while for the spiral compactum and the prolonged sine curve
we allow any finite Borel measure equivalent to \(H^1\). The argument for
the prolonged sine curve in fact applies to any reference measure assigning
positive measure to the limit interval \(I\).

\section{Root--contracting decompositions and Milnor attractors}\label{sec:root-contracting}

Let $(X,d)$ be a compact metric space and $A\subset X$ a nonempty closed
set.  Write
\[
X\setminus A \;=\; \bigsqcup_{\alpha\in \Theta} G_\alpha
\]
for the decomposition of the complement into its connected components.

Before introducing the abstract definition, we briefly explain the idea.
Given a closed set \(A\subset X\) and a connected component \(G\) of
\(X\setminus A\), we decompose \(G\) into open pieces, called
\emph{sectors}. Each sector has a chosen point \(a_i\in A\), called its
\emph{root}, which belongs to the closure of the sector. The closure of a
sector is allowed to meet \(A\) in more than one point. What matters is that
the dynamics inside the sector has all its \(\omega\)-limit sets contained in
\(A\), and that the chosen root \(a_i\) belongs to the \(\omega\)-limit set of
every interior point of the sector. The roots of the sectors are required to
be dense in the boundary of the component \(G\).

Near \(A\), the local maps move points by arbitrarily small amounts. This
ensures that the piecewise-defined global map is continuous on \(\cl(G)\) and
fixes \(A\) pointwise.

\begin{defn}[Root--contracting decomposition]\label{def:sector-decomposition}
Let \(X\) be a compact metric space and let \(A\subset X\) be closed. Write
\[
X\setminus A = \bigsqcup_{\alpha\in\Theta} G_\alpha
\]
for the decomposition of the complement into its connected components.
A \emph{root--contracting decomposition of \(X\setminus A\) toward \(A\)}
consists of a family of nonempty open sets \(\{U_{\alpha,i}\}_{\alpha\in\Theta,\,i\in I_\alpha}\),
with \(U_{\alpha,i}\subset G_\alpha\), points \(a_{\alpha,i}\in A\), and
continuous maps \(h_{\alpha,i}:\cl(U_{\alpha,i})\to\cl(U_{\alpha,i})\)
satisfying the following conditions.

\smallskip
\noindent\textbf{(1) Sector structure and fixed boundary.}
\begin{itemize}
    \item For each fixed \(\alpha\), the sets \(\{U_{\alpha,i}\}_{i\in I_\alpha}\)
    are pairwise disjoint and their closures cover \(G_\alpha\):
    \[
       G_\alpha \subset \bigcup_{i\in I_\alpha} \cl(U_{\alpha,i}).
    \]
    \item The chosen root belongs to the boundary impression of the sector:
    \[
       a_{\alpha,i}\in \cl(U_{\alpha,i})\cap A .
    \]
    \item For each \(\alpha\), the set \(\{a_{\alpha,i}\}_{i\in I_\alpha}\) is
    dense in \(\partial G_\alpha\).
    \item Each map \(h_{\alpha,i}\) fixes the boundary of \(U_{\alpha,i}\)
    pointwise:
    \[
       h_{\alpha,i}(x)=x \qquad\text{for all } x\in \cl(U_{\alpha,i})\setminus U_{\alpha,i}.
    \]
\end{itemize}

\smallskip
\noindent\textbf{(2) Rooted \(\omega\)-limit condition.}
For every \((\alpha,i)\) and every \(x\in U_{\alpha,i}\),
\[
        a_{\alpha,i}\in \omega_{h_{\alpha,i}}(x)\subset A .
\]

\smallskip
\noindent\textbf{(3) Slowing down near \(A\), uniformly over all components.}
For every \(b\in A\) and every \(\varepsilon>0\), there exists a neighbourhood
\(V\subset X\) of \(b\) such that, for \emph{every} \(\alpha\in\Theta\),
every \(i\in I_\alpha\), and every \(x\in U_{\alpha,i}\cap V\),
\[
d(h_{\alpha,i}(x),x)<\varepsilon .
\]
Note that \(V\) is required to depend only on \(b\) and \(\varepsilon\), not
on \(\alpha\): a single neighbourhood must control the displacement of the
local maps across every component of \(X\setminus A\) that meets it.

\smallskip
\noindent\textbf{(4) Local finiteness away from \(A\).}
For every \(x\in X\setminus A\), there exists a neighbourhood \(W\subset X\)
of \(x\) which meets only finitely many sets \(\cl(U_{\alpha,i})\).

\smallskip
If, in addition, each \(h_{\alpha,i}\) is a homeomorphism of
\(\cl(U_{\alpha,i})\) onto itself and the inverse maps also slow down near
\(A\) uniformly over all components — that is, for every \(b\in A\) and every
\(\varepsilon>0\), there exists a neighbourhood \(V\subset X\) of \(b\) such
that, for every \(\alpha\in\Theta\), every \(i\in I_\alpha\), and every
\(x\in U_{\alpha,i}\cap V\),
\[
d(h_{\alpha,i}^{-1}(x),x)<\varepsilon ,
\]
then the decomposition is called a \emph{homeomorphic root--contracting
decomposition (of \(X\setminus A\) toward \(A\))}.
\end{defn}

The following theorem gives the full-basin conclusion under a null-boundary
hypothesis. In the subsequent corollary, we assume in addition that the reference
measure assigns positive measure to every nonempty open set, which is used to
prove minimality of the attractor. This applies to the reference measures used
below, such as Lebesgue measure on compact manifolds and \(H^1\) on finite-length
one-dimensional continua.

\begin{thm}\label{thm:root-decomp-implies-attractor}
Let \(X\) be a compact locally connected metric space and let
\(A\subset X\) be closed. Let \(\{G_\alpha\}_{\alpha\in\Theta}\) be the
connected components of \(X\setminus A\). Assume that \(X\setminus A\) admits a root--contracting decomposition
toward \(A\) in the sense of Definition~\ref{def:sector-decomposition},
with associated data
\[
   \mathcal{U}
   =
   \{(U_{\alpha,i},a_{\alpha,i},h_{\alpha,i})\}_{\alpha\in\Theta,\ i\in I_\alpha}.
\]
(In practice this data is often built separately on each component
\(G_\alpha\); doing so gives a valid decomposition of \(X\setminus A\)
exactly when condition~(3) — the neighbourhood \(V\) in the slowing-down
axiom — can be taken independent of \(\alpha\).) Define a map \(f:X\to X\) by
\[
f(x)=
\begin{cases}
     h_{\alpha,i}(x), &
     x\in \cl(U_{\alpha,i}) \text{ for some }(\alpha,i),\\[2pt]
     x, & x\in A.
\end{cases}
\]
Then:
\begin{enumerate}
\item[\textnormal{(1)}] \textbf{Well-definedness and continuity.}
The map \(f\) is well defined and continuous on \(X\). Moreover, if the
decomposition is homeomorphic in the sense of
Definition~\ref{def:sector-decomposition}, then \(f\) is a homeomorphism of \(X\).
\item[\textnormal{(2)}] \textbf{Dynamics inside sectors.}
If \(x\in U_{\alpha,i}\) for some \((\alpha,i)\), then
\[
        \omega_f(x)=\omega_{h_{\alpha,i}}(x)\subset A
\]
and
\[
        a_{\alpha,i}\in \omega_f(x).
\]
If \(x\in \cl(U_{\alpha,i})\setminus U_{\alpha,i}\), then \(f(x)=x\).
\item[\textnormal{(3)}] \textbf{Full measure basin.}
Let \(\mu\) be any Borel probability measure on \(X\). If
\[
     \mu\!\left((\cl(U_{\alpha,i})\setminus U_{\alpha,i})\setminus A\right)=0
     \qquad\text{for all }(\alpha,i),
\]
then
\[
     \mu(B(A))=1.
\]
\end{enumerate}
\end{thm}

\begin{proof}
Each \(G_\alpha\) carries a root--contracting decomposition
\[
\{(U_{\alpha,i},a_{\alpha,i},h_{\alpha,i})\}_{i\in I_\alpha}
\]
satisfying the axioms of Definition~\ref{def:sector-decomposition}.

\medskip
\noindent\textbf{Step~1: Well--definedness.}
If \(x\in A\), then \(f(x)=x\). If \(x\in X\setminus A\), then \(x\) belongs to
a unique connected component \(G_\alpha\) of \(X\setminus A\). By the sector
covering axiom,
\[
G_\alpha\subset \bigcup_{i\in I_\alpha}\cl(U_{\alpha,i}),
\]
so \(x\in\cl(U_{\alpha,i})\) for at least one \(i\in I_\alpha\).

Suppose that
\[
x\in \cl(U_{\alpha,i})\cap\cl(U_{\beta,j}).
\]
If \(\alpha\neq\beta\), then the components \(G_\alpha\) and \(G_\beta\)
are open and closed in \(X\setminus A\). Hence
\(
\cl_X(G_\alpha)\cap(X\setminus A)=G_\alpha,
\)
so the closures of distinct components can meet only in \(A\). Therefore
\(x\in A\), and all
definitions give \(f(x)=x\). If \(\alpha=\beta\) and \(i\neq j\), then
\(x\notin U_{\alpha,i}\cup U_{\alpha,j}\), because the sectors are pairwise
disjoint. Hence
\[
x\in \cl(U_{\alpha,i})\setminus U_{\alpha,i}
\quad\text{and}\quad
x\in \cl(U_{\alpha,j})\setminus U_{\alpha,j}.
\]
Since each local map fixes the boundary of its sector pointwise, we get
\[
h_{\alpha,i}(x)=x=h_{\alpha,j}(x).
\]
Therefore the definition of \(f(x)\) is independent of the chosen sector. Thus
\(f\) is well defined.

\medskip
\noindent\textbf{Step 2: Continuity on \(X\setminus A\).}
Fix \(x\in X\setminus A\). By local finiteness, there exists a
neighbourhood \(W\subset X\setminus A\) of \(x\) meeting only finitely many
sector closures
\[
\cl(U_{\alpha_1,i_1}),\dots,\cl(U_{\alpha_m,i_m}).
\]
These closed sets cover \(W\), and on each
\(W\cap\cl(U_{\alpha_r,i_r})\) the map \(f\) agrees with the continuous
map \(h_{\alpha_r,i_r}\). Since the definitions agree on overlaps, the
finite gluing lemma shows that \(f|_W\) is continuous.

\medskip
\noindent\textbf{Step~3: Continuity at points of \(A\).}
Fix \(b\in A\) and let \(x_n\to b\) in \(X\). We show that \(f(x_n)\to b=f(b)\).

Let \(\varepsilon>0\). By condition~(3) of
Definition~\ref{def:sector-decomposition} — slowing down near \(A\),
\emph{uniformly over all components} \(G_\alpha\) — there exists a single
neighbourhood \(V\) of \(b\) such that, for \emph{every} \(\alpha\in\Theta\),
every \(i\in I_\alpha\), and every \(x\in U_{\alpha,i}\cap V\),
\[
d(h_{\alpha,i}(x),x)<\varepsilon/2.
\]
Crucially, \(V\) depends only on \(b\) and \(\varepsilon\), not on \(\alpha\);
this is precisely what condition~(3) provides, and it is what allows the
argument below to treat sectors from different components uniformly.

Shrinking \(V\) if necessary, we may also assume that
\[
d(y,b)<\varepsilon/2
\qquad\text{for all }y\in V.
\]
For all sufficiently large \(n\), we have \(x_n\in V\). We consider the three
possible locations of \(x_n\).

If \(x_n\in A\), then \(f(x_n)=x_n\), and hence
\[
d(f(x_n),b)=d(x_n,b)<\varepsilon/2.
\]

If \(x_n\in U_{\alpha,i}\) for some \((\alpha,i)\), then, since \(x_n\in
U_{\alpha,i}\cap V\), the displacement bound above applies directly:
\[
d(f(x_n),x_n)
=
d(h_{\alpha,i}(x_n),x_n)
<\varepsilon/2,
\]
and consequently
\[
d(f(x_n),b)
\le d(f(x_n),x_n)+d(x_n,b)
<\varepsilon/2+\varepsilon/2=\varepsilon.
\]

Finally, if \(x_n\in\cl(U_{\alpha,i})\setminus U_{\alpha,i}\) for some
\((\alpha,i)\), then \(f(x_n)=x_n\) by the boundary-fixing axiom, and again
\[
d(f(x_n),b)=d(x_n,b)<\varepsilon/2.
\]

In every case \(d(f(x_n),b)<\varepsilon\) for all sufficiently large \(n\).
Thus \(f(x_n)\to b\). Since \(b\in A\) was arbitrary, \(f\) is continuous at
every point of \(A\). Together with Step~2, this proves that \(f\) is
continuous on \(X\).

\smallskip\noindent
\textbf{Step 4: The homeomorphic case.}
Assume that the decomposition is homeomorphic. Define
\[
g(x)=
\begin{cases}
h_{\alpha,i}^{-1}(x),&x\in\cl(U_{\alpha,i}),\\
x,&x\in A.
\end{cases}
\]
Since each inverse map fixes the boundary of its sector pointwise, the same
argument as in Step~1 shows that \(g\) is well defined. The inverse maps satisfy the same local-finiteness and uniform
slowing-down conditions as the forward maps, so the arguments of
Steps~2 and~3 show that \(g\) is continuous. Finally,
\[
g\circ f=f\circ g=\operatorname{id}_X,
\]
since both identities hold on each sector closure and on \(A\). Thus
\(g=f^{-1}\), and \(f\) is a homeomorphism.

\medskip
\noindent\textbf{Step~5: Dynamics inside sectors.}
Let \(x\in X\setminus A\). Then \(x\in G_\alpha\) for a unique \(\alpha\), and
by the sector covering axiom (condition~(1) of
Definition~\ref{def:sector-decomposition}) \(x\in\cl(U_{\alpha,i})\) for some
\(i\in I_\alpha\).

If \(x\in U_{\alpha,i}\), then, since \(f\) agrees with \(h_{\alpha,i}\) on
\(\cl(U_{\alpha,i})\) and \(h_{\alpha,i}\) maps \(\cl(U_{\alpha,i})\) into
itself,
\[
f^n(x)=h_{\alpha,i}^{\,n}(x)
\qquad\text{for all }n\ge 0.
\]
Hence
\[
\omega_f(x)=\omega_{h_{\alpha,i}}(x).
\]
By the rooted \(\omega\)-limit axiom (condition~(2)),
\[
a_{\alpha,i}\in\omega_f(x)\subset A.
\]

If instead \(x\in\cl(U_{\alpha,i})\setminus U_{\alpha,i}\), then \(f(x)=x\) by
the boundary-fixing axiom. Thus \(x\) is a fixed point of \(f\); in
particular \(\omega_f(x)=\{x\}\), which need not lie in \(A\). These boundary
points are dealt with separately, and shown to be \(\mu\)-negligible, in
Step~6.

\medskip
\noindent\textbf{Step~6: The basin has full measure.}
Let \(\mu\) be a Borel probability measure on \(X\) satisfying
\[
\mu\!\left((\cl(U_{\alpha,i})\setminus U_{\alpha,i})\setminus A\right)=0
\qquad\text{for all }(\alpha,i).
\]

Every point of \(A\) lies in \(B(A)\), since \(f|_A=\operatorname{id}\), so
\(\omega_f(x)=\{x\}\subset A\) for \(x\in A\). Every point in
\[
\bigcup_{\alpha,i} U_{\alpha,i}
\]
belongs to \(B(A)\) by Step~5, since its \(\omega\)-limit set is contained in
\(A\).

We claim that the complement of these two sets,
\[
X\setminus\Big(A\cup\bigcup_{\alpha,i}U_{\alpha,i}\Big),
\]
is contained in
\[
\bigcup_{\alpha,i}
\left((\cl(U_{\alpha,i})\setminus U_{\alpha,i})\setminus A\right).
\]
Indeed, let \(x\) belong to this complement. Since \(x\notin A\), we have
\(x\in X\setminus A\), so \(x\in G_\alpha\) for some \(\alpha\), and by the
sector covering axiom \(x\in\cl(U_{\alpha,i})\) for some \(i\); since
\(x\notin\bigcup_{\alpha,i}U_{\alpha,i}\), in fact
\(x\in\cl(U_{\alpha,i})\setminus U_{\alpha,i}\). As also \(x\notin A\), this
gives \(x\in(\cl(U_{\alpha,i})\setminus U_{\alpha,i})\setminus A\), proving
the claim.

Since \(X\) is second countable and the sectors \(U_{\alpha,i}\) are pairwise
disjoint nonempty open sets, the index set \(\{(\alpha,i)\}\) is countable.
Hence the union on the right is a countable union of \(\mu\)-null sets, and
therefore itself has \(\mu\)-measure zero. Consequently,
\[
\mu\Big(X\setminus\big(A\cup\bigcup_{\alpha,i}U_{\alpha,i}\big)\Big)=0,
\]
and since \(A\cup\bigcup_{\alpha,i}U_{\alpha,i}\subset B(A)\), this gives
\[
\mu(X\setminus B(A))=0,
\qquad\text{so}\qquad
\mu(B(A))=1.
\]
\end{proof}

\begin{cor}\label{cor:gsp-implies-maf}
Let \(X\) be a compact locally connected metric space equipped with a
Borel probability measure \(\mu\) that assigns positive measure to every
nonempty open set. Suppose
that for every nonempty closed set \(A\subset X\), the complement
\(X\setminus A\) admits a root--contracting decomposition toward \(A\) in
the sense of Definition~\ref{def:sector-decomposition} — in particular, if
this decomposition is constructed separately on each connected component
\(G_\alpha\) of \(X\setminus A\), the resulting neighbourhoods in
condition~(3) (slowing down near \(A\)) must be chosen independently of
\(\alpha\), as in the remark following that definition. Assume moreover that for each \(A\) one can choose such a decomposition
so that the sector boundaries outside \(A\),
\[
        (\cl(U_{\alpha,i})\setminus U_{\alpha,i})\setminus A,
\]
have \(\mu\)-measure zero and the set of all chosen roots
\[
        \{a_{\alpha,i}:\alpha\in\Theta,\ i\in I_\alpha\}
\]
is dense in \(\partial A\). Then \(X\) is an MAF--space.

Moreover, if these decompositions can always be chosen homeomorphic in the
sense of Definition~\ref{def:sector-decomposition}, then every nonempty
closed set \(A\subset X\) can be realized as the Milnor global attractor of a
homeomorphism of \(X\).
\end{cor}

\begin{proof}
Let \(A\subseteq X\) be a nonempty closed set. By hypothesis, \(X\setminus A\)
admits a root--contracting decomposition toward \(A\). Applying
Theorem~\ref{thm:root-decomp-implies-attractor} yields a continuous map
\(f_A:X\to X\) such that
\[
        \mu(B(A))=1.
\]
In particular,
\[
        \omega(x)\subseteq A
\]
for \(\mu\)-almost every \(x\in X\), and hence the likely limit set
\(\Lambda(f_A)\) is a closed subset of \(A\).

If, moreover, the decomposition can be chosen homeomorphic in the sense of
Definition~\ref{def:sector-decomposition}, then
Theorem~\ref{thm:root-decomp-implies-attractor} gives \(f_A\) as a
homeomorphism.

\medskip
\noindent\textbf{Minimality of \(A\).}
Let \(B\subsetneq A\) be closed. Choose \(a\in A\setminus B\), and let
\(U\subset X\setminus B\) be an open neighbourhood of \(a\) (such a
neighbourhood exists since \(B\) is closed and \(a\notin B\)). Since \(A\)
is closed, \(A=\operatorname{int}(A)\sqcup\partial A\), so exactly one of the
following two cases occurs.

If \(U\cap\operatorname{int}(A)\neq\varnothing\), choose a nonempty open set
\(W\subset U\cap\operatorname{int}(A)\). Since \(f_A\) fixes \(A\) pointwise,
\(\omega(x)=\{x\}\) for every \(x\in W\subset A\), and \(x\notin B\), so
\(\omega(x)\not\subset B\). As \(\mu(W)>0\) by hypothesis on \(\mu\), the set
\(B\) cannot have full basin.

Otherwise \(a\in\partial A\). Since the set of all chosen roots
\(\{a_{\alpha,i}:\alpha\in\Theta,\ i\in I_\alpha\}\) is dense in \(\partial
A\) and \(U\) is a neighbourhood of \(a\in\partial A\), there are
\(\alpha\in\Theta\) and \(i\in I_\alpha\) with \(a_{\alpha,i}\in U\). Since
\(a_{\alpha,i}\in\cl(U_{\alpha,i})\) and \(U\) is an open neighbourhood of
\(a_{\alpha,i}\), every neighbourhood of \(a_{\alpha,i}\) meets
\(U_{\alpha,i}\); hence \(U_{\alpha,i}\cap U\) is nonempty. It is also open,
so it has positive \(\mu\)-measure. For every \(x\in U_{\alpha,i}\cap U\),
the rooted \(\omega\)-limit condition gives \(a_{\alpha,i}\in\omega(x)\).
Since \(a_{\alpha,i}\in U\subset X\setminus B\), we have
\(a_{\alpha,i}\notin B\), so \(\omega(x)\not\subset B\). Thus \(B\) cannot
have full basin.

In either case, \(\mu(B(B))<1\). We have shown that no proper closed subset of \(A\) has full basin. If
\(B\subset X\) is any closed set with full basin, then \(A\cap B\) also
has full basin, since almost every \(\omega\)-limit set is contained in
both \(A\) and \(B\). Hence \(A\cap B=A\), so \(A\subset B\). Therefore
\(A\) is the smallest closed full-basin set, and
\(\Lambda(f_A)=A\). Since \(A\) was arbitrary, \(X\)
is an MAF--space.
\end{proof}

\section{One-dimensional manifolds}\label{sec:1D-manifolds}

Throughout this section, \(\lambda\) denotes the normalized Lebesgue
(length) measure induced by a fixed Riemannian metric on \(M\).

\begin{thm}[Milnor attractors on one-dimensional manifolds]\label{thm:1D}
Let \(M\) be a compact one-dimensional manifold, with or without boundary,
equipped with the measure \(\lambda\). Then every nonempty closed subset
\(A\subset M\) is the Milnor global attractor of a continuous self-map of
\(M\).

If \(M\) is connected, then the realizing map may be chosen to be a
homeomorphism.
\end{thm}

\emph{Idea of the proof.}
The connected components of $M\setminus A$ are intervals. If both endpoints of
such an interval belong to $A$, we split it into two sectors at its midpoint;
otherwise the whole interval forms a single sector. Each sector therefore has a
distinguished root belonging to $A$.

We then conjugate a fixed model homeomorphism of the open interval to every
sector. The model map moves every point monotonically toward the root while
slowing down near the boundary of the sector. Gluing these maps with the
identity on $A$ yields a homeomorphic root--contracting decomposition of
$M\setminus A$ toward $A$. The abstract realization theorem then implies that
$A$ is the Milnor global attractor.

The purpose of presenting the one-dimensional case first is that exactly the
same mechanism will be used in the next section for surfaces. There, intervals
are replaced by sectors obtained from a suitable decomposition of the
complement, while the model dynamics and the gluing argument remain the same in
spirit.

\begin{lemma}\label{lem:1d-accessible}
Let $M$ be a compact one-dimensional manifold and let
$A\subset M$ be closed. Write
$$
M\setminus A
=
\bigsqcup_{\alpha\in\Theta}G_\alpha,
$$
where the $G_\alpha$ are the connected components. Let
$$
\operatorname{AC}(A)
=
\bigcup_{\alpha\in\Theta}
\operatorname{AC}(G_\alpha),
$$
where $\operatorname{AC}(G_\alpha)$ denotes the set of points
$p\in\partial G_\alpha$ for which there exists an embedded (i.e.\ injective
continuous) arc
$$
\gamma:[0,1]\to M
$$
such that
$$
\gamma(0)=p
\qquad\text{and}\qquad
\gamma((0,1])\subset G_\alpha.
$$
Then $\operatorname{AC}(A)$ is dense in $\partial A$.
Consequently, $\partial A$ contains a countable dense subset whose
points are accessible from complementary components.
\end{lemma}

\begin{proof}
Let $a\in\partial A$ and let $U$ be a neighbourhood of $a$. Choose a chart \(\varphi:W\to I\), with \(\varphi(a)=0\) and \(W\subset U\),
where \(I=(-\delta,\delta)\) if \(a\notin\partial M\), and
\(I=[0,\delta)\) if \(a\in\partial M\). Choose
\(x\in W\setminus A\). After reversing the coordinate in the first case
if necessary, write \(t=\varphi(x)>0\).

The set $K=\varphi(A\cap W)\cap[0,t]$ is compact (as $A\cap W$ is closed in
$W$ and $\varphi$ is a homeomorphism) and nonempty (it contains
$\varphi(a)=0$). Let $s_0=\max K$ and $p=\varphi^{-1}(s_0)\in A$.

By maximality of $s_0$, $(s_0,t]\cap K=\varnothing$, so
$\varphi^{-1}((s_0,t])\subset W\setminus A\subset M\setminus A$. This set is
connected and contains $x$, hence lies in the component $G_\alpha$ of
$M\setminus A$ containing $x$. Since $p=\lim_{r\to0^+}\varphi^{-1}(s_0+r)$
with $\varphi^{-1}(s_0+r)\in G_\alpha$, and $p\in A$, we get
$p\in\operatorname{cl}(G_\alpha)\setminus G_\alpha=\partial G_\alpha$.

The arc $\gamma(r)=\varphi^{-1}(s_0+r(t-s_0))$, $r\in[0,1]$, is embedded,
with $\gamma(0)=p$ and $\gamma((0,1])=\varphi^{-1}((s_0,t])\subset G_\alpha$.
Thus $p\in\operatorname{AC}(G_\alpha)\subset\operatorname{AC}(A)$, and since
$s_0\in[0,t)\subset(-\delta,\delta)$, we have $p\in W\subset U$, so
$U\cap\operatorname{AC}(A)\neq\varnothing$.

Since $U$ was an arbitrary neighbourhood of $a\in\partial A$,
$\operatorname{AC}(A)$ is dense in $\partial A$.
\end{proof}

\begin{proposition}\label{prop:1d-root-decomposition}
Let $M$ be a compact connected one-dimensional manifold and let
$A\subsetneq M$ be closed and nonempty. Then there exist a nonempty at most countable index set \(N\), a family
of pairwise disjoint nonempty open intervals \(\{G_i\}_{i\in N}\)
that is locally finite in \(M\setminus A\), roots \(g_i\in\partial A\),
and homeomorphisms \(f_i:\cl(G_i)\to\cl(G_i)\) such that:
\begin{enumerate}
\item the closures of the sectors cover $M\setminus A$, and the set
$\{g_i:i\in N\}$ is dense in $\partial A$;
\item $f_i$ fixes $\cl(G_i)\setminus G_i$ pointwise, and $f_i^{\,n}(x)\to g_i$
for every $x\in G_i$;
\item the families $\{f_i\}$ and $\{f_i^{-1}\}$ slow down near $A$
uniformly in $i$;
\item $\lambda\bigl((\cl(G_i)\setminus G_i)\cap(M\setminus A)\bigr)=0$ for
every $i$.
\end{enumerate}
Consequently, $\{(G_i,g_i,f_i)\}_{i\in N}$ is a homeomorphic
root--contracting decomposition of $M\setminus A$ toward $A$.
\end{proposition}

\begin{proof}
Every connected component $I$ of $M\setminus A$ is, as a subset of $M$,
either an open interval whose two topological ends both belong to $A$
(necessarily so, by compactness and maximality of $I$ as a component,
unless $I$ meets $\partial M$), or — when $M$ has nonempty boundary and $I$
meets $\partial M$ — a half-open interval $I=J\sqcup\{p_I\}$, where
$p_I\in\partial M\setminus A$ is an actual point of $I$ and $J=I\setminus
\{p_I\}$ is open, with the other topological end of $J$ belonging to $A$.

In the first case, if the two ends of $I$ correspond to distinct points of
$A$, or to the same point of $A$, choose $m_I\in I$ and split $I$ into two
sectors at $m_I$, taking these point(s) as roots (both sectors get the same
root in the latter case); the point $m_I$ lies in neither open sector but
in the closure of both, where it will be fixed pointwise by both local
maps. In the second case, take the sector to be $G_i:=J=I\setminus\{p_I\}$,
rooted at the unique point $g_i\in A$ that is the other topological end of
$J$; the point $p_I$ then plays, for this single sector, exactly the role
that $m_I$ plays for a split pair — it lies in $\cl(G_i)\setminus G_i$, not
in $G_i$ itself, so it is exempt from the rooted $\omega$-limit condition.
The resulting sectors are pairwise disjoint open intervals, their closures
cover $M\setminus A$, and the family is locally finite in $M\setminus A$.
Every point of $\operatorname{AC}(A)$ occurs as a root, so by
Lemma~\ref{lem:1d-accessible} the roots are dense in $\partial A$.

Let $B(t)=\exp(-1/(1-t^2))$ for $|t|<1$ and $B(t)=0$ for $|t|\geq1$, and set
$\psi(x)=B(2x-1)$. Choose $\varepsilon>0$ small enough that
$f_I(x)=x+\varepsilon\psi(x)$ is strictly increasing on $[0,1]$; then $f_I$
is a homeomorphism of $[0,1]$ fixing $0$ and $1$, with $f_I^{\,n}(x)\to1$
for every $x\in(0,1)$. For each $i$, choose a homeomorphism
$\phi_i:[0,1]\to\cl(G_i)$ with $\phi_i(1)=g_i$, and set
$f_i=\phi_i\circ f_I\circ\phi_i^{-1}$. Then $f_i$ fixes $\cl(G_i)\setminus
G_i$ pointwise — this set consists of the root $g_i$ together with one
further point, $m_I$ or $p_I$ according to which case $G_i$ arose from,
and $f_I$ fixes both endpoints $0,1$ of $[0,1]$ — and $f_i^{\,n}(x)\to g_i$
for every $x\in G_i$.

For the slowing-down property, fix $b\in A$ and $\varepsilon_0>0$. Since $M$
has finite length, only finitely many complementary intervals
$G_{i_1},\dots,G_{i_k}$ have diameter at least $\varepsilon_0$. Among these,
at most two contain $b$ in their boundary, since $b$ has at most two adjacent
components. For each such $G_{i_j}$, continuity of $f_{i_j}^{\pm1}$ at the
fixed point $b$ gives a neighbourhood $V_j$ of $b$ on which
$d(f_{i_j}^{\pm1}(x),x)<\varepsilon_0$. The remaining large intervals have
closures disjoint from $b$, so, since there are only finitely many of them,
there exists a neighbourhood $V_0$ of $b$ disjoint from all their closures.
Set
\[
V=V_0\cap\bigcap_j V_j,
\]
which is again a neighbourhood of $b$. Every sector meeting $V$ then either
has diameter less than $\varepsilon_0$ (so
$d(f_i^{\pm1}(x),x)\le\operatorname{diam}(G_i)<\varepsilon_0$ directly) or is
one of the intervals $G_{i_j}$ containing $b$ in its boundary, where the
bound holds on $V\subset V_j$ by construction. Since $\varepsilon_0$ and $b$
were arbitrary, $\{f_i\}$ and $\{f_i^{-1}\}$ slow down near $A$, uniformly
in $i$.

Finally, each sector boundary $\cl(G_i)\setminus G_i$ consists of the root
$g_i$ together with exactly one further point — $m_I$ in the split case,
$p_I$ in the boundary case — and that point does not lie in $A$. Hence
$\lambda\bigl((\cl(G_i)\setminus G_i)\cap(M\setminus A)\bigr)=0$.
\end{proof}
\begin{proof}[Proof of Theorem~\ref{thm:1D}]
Let \(A\subset M\) be nonempty and closed. If \(A=M\), the identity map
realizes \(A\), so assume \(A\neq M\). Let \(M_1,\dots,M_r\) be the
connected components of \(M\), and set \(A_\alpha=A\cap M_\alpha\).

If \(A_\alpha=M_\alpha\), let
\(f_\alpha=\operatorname{id}_{M_\alpha}\). If
\(\varnothing\neq A_\alpha\subsetneq M_\alpha\), apply
Proposition~\ref{prop:1d-root-decomposition} and
Corollary~\ref{cor:gsp-implies-maf} to obtain a homeomorphism
\(f_\alpha:M_\alpha\to M_\alpha\) with Milnor global attractor
\(A_\alpha\).

Choose \(a\in A\), and on each component disjoint from \(A\) set
\(f_\alpha\equiv a\). Since the components are clopen, these maps combine
into a continuous map \(f:M\to M\), which is a homeomorphism whenever
\(A\) meets every component. Every point in a component disjoint from
\(A\) satisfies \(\omega(x)=\{a\}\subset A\), while \(A_\alpha\) has full
basin in every component meeting \(A\). Hence \(\lambda(B(A))=1\).

For minimality, let \(B\subsetneq A\) be closed. Then
\(B\cap M_\alpha\subsetneq A_\alpha\) for some component \(M_\alpha\)
meeting \(A\). Since \(A_\alpha\) is the Milnor global attractor of
\(f_\alpha\), the set \(B\cap M_\alpha\) does not have full basin in
\(M_\alpha\). Hence \(B\) does not have full basin in \(M\).

Now let \(K\subset M\) be any closed set with full basin. Then
\(A\cap K\) also has full basin, since almost every \(\omega\)-limit set
is contained in both \(A\) and \(K\). By the preceding argument,
\(A\cap K\) cannot be a proper subset of \(A\). Thus \(A\subset K\), so
\(A\) is the smallest closed set with full basin and therefore
\(\Lambda(f)=A\).
\end{proof}
\begin{rem}
    In dimension one we allow manifolds with boundary, since the interval
endpoints require only minor modifications of the construction. In
dimension two we restrict to closed manifolds. Treating boundary
components requires additional control of the triangulation and sector
dynamics near the boundary. We expect the realization theorem to remain
valid for compact surfaces with boundary, but do not pursue this extension
here.
\end{rem}
\section{Two-dimensional manifolds}\label{sec:smooth_compact_manifolds}
We now turn to compact two-dimensional manifolds, where $\lambda$ denotes
the Lebesgue measure induced by a Riemannian metric.

\begin{thm}[Milnor attractors on two-dimensional manifolds]\label{thm:2D}
Let $S$ be a smooth closed connected \(2\)-dimensional manifold equipped
with Lebesgue measure \(\lambda\). Then for every nonempty closed subset
\(A\subset S\), there exists a homeomorphism \(f:S\to S\) for which \(A\) is
the Milnor global attractor with respect to \(\lambda\). In particular,
every smooth compact connected \(2\)-manifold is an MAF--space.

More generally, if $S$ is not connected, then every nonempty closed subset
\(A\subset S\) is the Milnor global attractor of some continuous map
\(f:S\to S\), and \(f\) can be chosen to be a homeomorphism whenever \(A\)
meets every connected component of $S$.
\end{thm}

\noindent\emph{Idea of the proof.}
Let $G$ be a connected component of $S\setminus A$. We choose a countable
dense set of accessible points of $\partial G$ and construct pairwise
disjoint end-cuts landing at these points. An adapted shrinking
triangulation of $G$ and Whitehead's lemma then produce a spanning forest
of special trees in the dual graph. Propagating colors through these trees
decomposes $G$ into sectors, each homeomorphic to an open disk and equipped
with a proper ray accumulating at a prescribed root in $\partial G$. The
shrinking triangulation ensures local finiteness, zero-measure sector
boundaries outside $A$, and that consecutive layers lie in the same or
adjacent triangles.

A model homeomorphism of the open triangle is transferred to each sector.
Its orbits have $\omega$-limit sets contained in $A$ and containing the
distinguished root, while the shrinking triangulation ensures that the
sector maps slow down near $A$.
Proposition~\ref{prop:verify-axioms} shows that the resulting
componentwise constructions satisfy the axioms of a homeomorphic
root--contracting decomposition of $S\setminus A$ toward $A$.
Theorem~\ref{thm:root-decomp-implies-attractor} and
Corollary~\ref{cor:gsp-implies-maf} then imply that $A$ is the Milnor
global attractor of a homeomorphism. If $S$ is disconnected, the
construction is carried out on each component meeting $A$, while components
disjoint from $A$ are mapped constantly into $A$, yielding a continuous
realizing map.

\begin{lemma}[Gluing disks]\label{lem:gluing}
Let \(D_1,D_2\) be closed topological disks such that
\(D_1\cap D_2\) is a nondegenerate boundary arc of each.
Then \(D_1\cup D_2\) is a closed topological disk.

Let \(E_1,E_2\) be another such pair, and let
\(h:D_1\to E_1\) be a homeomorphism with
\[
h(D_1\cap D_2)=E_1\cap E_2.
\]
Then \(h\) extends to a homeomorphism
\[
H:D_1\cup D_2\to E_1\cup E_2.
\]
The extension may be chosen to match a prescribed boundary arc
of \(D_2\), disjoint from \(D_1\cap D_2\), with a prescribed
boundary arc of \(E_2\), disjoint from \(E_1\cap E_2\).
It may simultaneously match a prescribed properly embedded arc
joining the relative interiors of the attaching and prescribed
boundary arcs in \(D_2\) with such an arc in \(E_2\), provided
their endpoints on the attaching arcs correspond under \(h\).
\end{lemma}

\begin{proof}
The first assertion follows by identifying \(D_1\) and
\(D_2\) with two half-disks whose common boundary diameter
is identified.

Let
\[
I=D_1\cap D_2,
\qquad
J=E_1\cap E_2.
\]
Extend \(h|_I\) to a homeomorphism
\(
\widetilde h:\partial D_2\to\partial E_2
\)
with the prescribed boundary-arc and endpoint
correspondences. Cutting along the prescribed properly
embedded arc, if present, and applying the Schoenflies
extension theorem gives a homeomorphism
\(h_2:D_2\to E_2\) with the required correspondences.
Since \(h=h_2\) on \(I\), the two maps glue to the required
homeomorphism \(H\).
\end{proof}

The following lemma is an adaptation of a well-known result on accessible points of plane domains due to Newman \cite{Newman}.

\begin{lemma}\label{lem:2d-accessible}
Let $S$ be a smooth compact two-dimensional manifold and let $A\subset S$ be
closed. Write $S\setminus A=\bigsqcup_{\alpha\in\Theta}G_\alpha$ for the
decomposition into connected components, and let
$\operatorname{AC}(A)=\bigcup_\alpha\operatorname{AC}(G_\alpha)$, where
$\operatorname{AC}(G_\alpha)$ is the set of $p\in\partial G_\alpha$ for which
some simple arc $L\subset S$ has $p\in L$ and $L\setminus\{p\}\subset
G_\alpha$. Then $\operatorname{AC}(A)$ is dense in $\partial A$;
consequently $\partial A$ contains a countable dense set whose points are
accessible from components of $S\setminus A$.
\end{lemma}

\begin{proof}
Since $S$ is locally connected, every component $G_\alpha$ of the open set
$S\setminus A$ is open; being also closed in $S\setminus A$, it satisfies
$\partial G_\alpha\subset A$, so $\operatorname{AC}(A)\subset\partial A$.

Let $a\in\partial A$ and $\varepsilon>0$, and choose a connected open
neighbourhood $a\in V\subset B(a,\varepsilon)$. Since $a\in\partial A$,
there is $x\in V\setminus A$; let $G_\alpha$ be the component of
$S\setminus A$ containing $x$. As connected open subsets of a
two-dimensional manifold are arcwise connected, there is a simple arc
$\gamma:[0,1]\to V$ with $\gamma(0)=a$, $\gamma(1)=x$.

Let $E=\{t\in[0,1]:\gamma([t,1])\subset G_\alpha\}$ and $t_0=\inf E$. Since
$G_\alpha$ is open and $\gamma(1)=x\in G_\alpha$, $E\neq\varnothing$; for
every $t>t_0$ there is $s\in E$ with $s<t$, so $\gamma(t)\in G_\alpha$,
giving $\gamma((t_0,1])\subset G_\alpha$. If $\gamma(t_0)\in G_\alpha$, then
openness of $G_\alpha$ and continuity of $\gamma$ would give $t<t_0$ with
$\gamma([t,t_0])\subset G_\alpha$, hence $t\in E$, contradicting $t_0=\inf
E$. So $\gamma(t_0)\in\cl(G_\alpha)\setminus G_\alpha=\partial G_\alpha$, and
the subarc $\gamma([t_0,1])$ exhibits $\gamma(t_0)$ as accessible from
$G_\alpha$, i.e. $\gamma(t_0)\in\operatorname{AC}(A)$. Since
$\gamma(t_0)\in V\subset B(a,\varepsilon)$, this shows $\operatorname{AC}(A)$
is dense in $\partial A$.

Finally, let $\{U_n:n\in\mathbb N\}$ be a countable basis for $\partial A$
(second countable, as a subset of $S$), and for each $n$ with
$U_n\cap\operatorname{AC}(A)\neq\varnothing$ choose $a_n\in
U_n\cap\operatorname{AC}(A)$. Then $K=\{a_n:U_n\cap\operatorname{AC}(A)
\neq\varnothing\}$ is countable, dense in $\partial A$ by density of
$\operatorname{AC}(A)$, and every point of $K$ is accessible from some
component of $S\setminus A$.
\end{proof}

A standard tool in the study of two-dimensional manifolds is the use of
triangulations, which allow one to replace the manifold by a combinatorial
structure of simplices without loss of topological information.

\begin{defn}
Let $X$ be a two-dimensional manifold. A \emph{triangulation} of $X$ is a
countable, locally finite collection of closed topological triangles
covering $X$, such that distinct triangles meet only in a full edge or a
vertex. A \emph{topological triangle} in $X$ is a continuous embedding
$\mathbb{B}\hookrightarrow X$, where $\mathbb{B}\subset\mathbb{R}^2$ is a
closed triangle; for each topological triangle $\Delta\subset X$, we write
$\phi_\Delta:\Delta\to\mathbb{B}$ for a homeomorphism taking vertices to
vertices. A triangulation is \emph{smooth} if each defining embedding
$\mathbb{B}\hookrightarrow X$ is a smooth embedding, in which case each
$\phi_\Delta$ may be taken to be a diffeomorphism.
\end{defn}

By a theorem of Radó (1925), every two-dimensional manifold is
triangulable; see also Moise (1977). The triangulation can moreover be
chosen smooth, see Munkres~\cite[Thm.~10.6]{Munkres}.

For a triangulation \(T\) of \(X\), the \emph{dual graph}
\(D(T)\) has vertex set
\(
V(D(T))=\{\Delta:\Delta\in T\},
\)
with two vertices joined by an edge whenever the corresponding
triangles share a side. To realize \(D(T)\) in \(X\), choose
an interior point of each shared side and represent the vertex
corresponding to \(\Delta\) by the inverse image under
\(\phi_\Delta\) of the barycentre of \(\mathbb B\).
Join this point to the chosen points on the sides using
straight segments in \(\mathbb B\), pulled back by
\(\phi_\Delta^{-1}\). The two segments associated with a
shared side form the corresponding dual edge.

A \emph{tree} is a connected graph in which any two vertices
are joined by a unique simple path. A \emph{ray} is a
one-way infinite simple path \(V_0,V_1,V_2,\ldots\), with
\(V_iV_{i+1}\) an edge for every \(i\geq0\). For a subgraph
\(H\subset D(T)\), write \(|H|\) for its geometric realization.

Let \(G\) be an open subset of \(S\).
A \emph{proper topological ray} in \(G\) is an embedding
\(\rho:[0,\infty)\to G\) that eventually leaves every
compact subset of \(G\). Its accumulation set in \(S\) is
\[
\operatorname{Acc}(\rho)
=\bigcap_{t\geq0}\cl_S\bigl(\rho([t,\infty))\bigr).
\]
For a graph ray \(R\) whose geometric realization is proper
in \(G\), we also write \(\operatorname{Acc}(|R|)\) for this
set, using any parametrization of \(|R|\) along the ray.

A \emph{special tree} is a locally finite tree of the form
\[
C=R\cup\bigcup_{\alpha\in I}F_\alpha,
\]
where \(R\) is a ray, called the \emph{distinguished ray},
and the \(F_\alpha\) are pairwise disjoint finite trees,
each meeting \(R\) in exactly one vertex.
A \emph{forest} is a graph whose components are trees.

Whitehead~\cite[Lemma~3.1]{Whitehead} proved that every connected
infinite locally finite graph admits a spanning forest whose
components are special trees. The following lemma adapts this
decomposition to the dual graph of a triangulation of \(G\),
with prescribed distinguished rays landing at accessible
boundary points.

\begin{lemma}[Whitehead forest]\label{lem:whitehead-forest}
Let \(S\) be a smooth compact two-dimensional manifold, let
\(A\subset S\) be closed, and let \(G\) be a component of
\(S\setminus A\) with \(\partial G\neq\varnothing\). Let
\(Q=\{q_n:n\in N\}\subset\operatorname{AC}(G)\), where \(N\) is a
nonempty finite or countable index set and the points \(q_n\) are
distinct. We take \(N=\{1,\dots,r\}\) for some \(r\geq1\), or \(N=\mathbb N\). Then there exist:
\begin{enumerate}
\item pairwise disjoint end-cuts $\gamma_n:[0,1]\to\cl(G)$ with
$\gamma_n([0,1))\subset G$ and $\gamma_n(1)=q_n$;
\item a countable locally finite smooth triangulation \(\mathcal T\) of \(G\)
such that
\[
\operatorname{diam}(\Delta)
<\tfrac12\operatorname{dist}(\Delta,A)
\]
for every \(\Delta\in\mathcal T\);
\item pairwise vertex-disjoint rays \(R_n\subset D(\mathcal T)\)
whose geometric realizations are proper in \(G\) and satisfy
\(\operatorname{Acc}(|R_n|)=\{q_n\}\), where accumulation sets
are taken in \(S\);
\item a spanning forest $F\subset D(\mathcal T)$ whose components are
special trees, such that every $R_n$ is the distinguished ray of a
component of $F$.
\end{enumerate}
In particular, if \(Q\) is dense in \(\partial G\), then the
distinguished rays of \(F\) include rays landing at a dense
subset of \(\partial G\).
\end{lemma}

\begin{proof}
First we choose the end-cuts inductively, with
\(\gamma_n([0,1])\subset B(q_n,2^{-n})\).
The union of the previously chosen end-cuts is compact and
does not contain \(q_n\), so a sufficiently short terminal
subarc of an end-cut landing at \(q_n\) avoids this union
and lies in \(B(q_n,2^{-n})\). Reparametrize this short subarc
as \(\gamma_n:[0,1]\to\cl(G)\).

The sets \(E_n=\gamma_n([0,1))\) are pairwise disjoint
and closed in \(G\). To see that they form a locally finite
family, fix \(x\in G\) and choose a closed ball
\(K=\overline B(x,\varepsilon)\subset G\). Since \(K\) is
compact,
\(
\delta=\operatorname{dist}(K,\partial G)>0.
\)
As \(E_n\subset B(q_n,2^{-n})\) and \(q_n\in\partial G\),
we have \(E_n\cap K=\varnothing\) whenever \(2^{-n}<\delta\).
Thus the neighbourhood \(B(x,\varepsilon)\) meets only
finitely many \(E_n\).

We now choose a triangulation separating the end-cuts.
By local finiteness, each \(x\in G\) has a neighbourhood
meeting only finitely many \(E_n\). Since these sets are
closed in \(G\) and pairwise disjoint, we may shrink this
neighbourhood to meet at most one of them. Choose
\(r_x>0\) so that \(B(x,r_x)\) lies in this neighbourhood
and \(5r_x<d(x,A)\).

Start with a countable locally finite smooth triangulation
of \(G\). For each triangle \(\Delta\), the sets
\(\phi_\Delta(\Delta\cap B(x,r_x))\) form an open cover
of the compact model triangle \(\mathbb B\). Hence there
is \(\delta_\Delta>0\) such that every subset of
\(\mathbb B\) of diameter less than \(\delta_\Delta\)
lies in one member of this cover. After finitely many
barycentric subdivisions, every resulting closed triangle
has diameter less than \(\delta_\Delta\). Pulling this
subdivision back by \(\phi_\Delta^{-1}\), we obtain a
finite subdivision of \(\Delta\), each triangle of which
is contained in some \(B(x,r_x)\).
On each original edge, retain all division points introduced
from either side, so that both incident triangles use the
same partition. In the image of each small triangle under
\(\phi_\Delta\), choose an interior point and join it by
straight segments to all boundary vertices, including the
additional division points. Pulling these subdivisions back
gives compatible subdivisions of the original triangles.
Every new triangle remains contained in one of the
previously chosen small triangles. Each original triangle has been
subdivided into finitely many triangles, so the resulting
triangulation \(\mathcal T\) is countable, locally finite,
and smooth.

Every \(\Delta\in\mathcal T\) lies in some \(B(x,r_x)\),
and therefore meets at most one end-cut and satisfies
\[
\operatorname{diam}(\Delta)
\le 2r_x
<\tfrac12\bigl(d(x,A)-r_x\bigr)
\le\tfrac12\operatorname{dist}(\Delta,A).
\]

For each \(n\), let \(H_n\) be the dual subgraph induced
by the triangles meeting \(E_n\). It is connected:
a compact subarc joining points in two such triangles
connects them through a finite chain of triangles.
Crossing an edge gives a dual edge; at a triangulation
vertex on the subarc, all incident triangles meet \(E_n\)
and are connected by the dual edges around that vertex.
Moreover, \(H_n\) is infinite, since a finite union of
triangles is compact in \(G\) and cannot contain \(E_n\).

Thus \(H_n\) is connected, infinite, and locally finite,
so it contains a ray \(R_n\). The rays are pairwise
vertex-disjoint because no triangle meets two end-cuts.
Their geometric realizations are proper in \(G\), by
local finiteness of \(\mathcal T\).

Write \(\Delta_{n,k}\) for the successive distinct triangles
along \(R_n\), and choose \(x_{n,k}\in\Delta_{n,k}\cap E_n\).
Every compact initial subarc of \(\gamma_n\) meets only
finitely many triangles, so \(x_{n,k}\to q_n\). Since
\[
\operatorname{diam}(\Delta_{n,k})
<\tfrac12\operatorname{dist}(\Delta_{n,k},A)
\le\tfrac12d(x_{n,k},q_n)\longrightarrow0,
\]
the tails of \(|R_n|\) converge to \(q_n\). Hence
\(\operatorname{Acc}(|R_n|)=\{q_n\}\).

Finally, extend \(\{R_n:n\in N\}\) to a maximal family
\(\mathscr R\) of pairwise vertex-disjoint rays in
\(D(\mathcal T)\), and let \(P\) be their combined vertex
set. Every component of the induced graph on
\(V(D(\mathcal T))\setminus P\) is finite, since an infinite
component would contain another ray. By connectedness of
\(D(\mathcal T)\), each such component has an edge to \(P\).
Choose a spanning tree in each component and attach it
to one ray by exactly one such edge.

Keeping only the ray edges and these attached trees gives
a spanning forest. Local finiteness ensures that only
finitely many trees attach at each ray vertex; merging
them gives a single finite tree there. Every component
is therefore a special tree, and each prescribed \(R_n\)
remains the distinguished ray of a distinct component.
The density assertion follows from
\(\operatorname{Acc}(|R_n|)=\{q_n\}\).
\end{proof}

A \emph{polygonal arc} in \(G\) is a finite concatenation of smooth
embedded arcs. A \emph{polygonal ray} is a proper embedding
\(\rho:[0,\infty)\to G\) whose restriction to every compact interval is
polygonal. A \emph{polygonal disk} is a compact subset \(D\subset G\)
homeomorphic to the closed unit disk and whose boundary is polygonal.

\begin{lemma}[Propagation through a triangle]
\label{lem:triangle-propagation}
Let \(\Delta\) be a triangle of a smooth triangulation, with a chosen side
\(e_{\mathrm{out}}\), and suppose that the other two sides are partitioned
into nondegenerate closed arcs \(I_1,\dots,I_r\) with pairwise disjoint
interiors. Then there exist nondegenerate closed arcs
\(J_1,\dots,J_r\subset e_{\mathrm{out}}\), with pairwise disjoint interiors
and union \(e_{\mathrm{out}}\), and polygonal disks
\(P_1,\dots,P_r\subset\Delta\), with pairwise disjoint interiors, such that
\[
\Delta=\bigcup_{j=1}^r P_j
\qquad\text{and}\qquad
P_j\cap\partial\Delta=I_j\cup J_j
\]
for every \(j\). Moreover, the order of the \(J_j\) along
\(e_{\mathrm{out}}\) agrees with the order of the \(I_j\) along the other
two sides.
\end{lemma}
\begin{proof}
Let $p_0,p_1,\dots,p_r$ be the common endpoints of consecutive arcs among
$I_1,\dots,I_r$, in order along the two sides of $\Delta$ other than
$e_{\mathrm{out}}$ (so $I_j$ runs from $p_{j-1}$ to $p_j$), with $p_0,p_r$
the endpoints of $e_{\mathrm{out}}$. Choose a partition $e_{\mathrm{out}}=
J_1\cup\cdots\cup J_r$ into nondegenerate closed arcs with pairwise disjoint
interiors, in this same order, and let $q_0,q_1,\dots,q_r$ be the resulting
division points along $e_{\mathrm{out}}$ (so $J_j$ runs from $q_{j-1}$ to
$q_j$, with $q_0=p_0$ and $q_r=p_r$ the endpoints of $e_{\mathrm{out}}$).

Fix a diffeomorphism $\phi_\Delta:\Delta\to B$ taking vertices to
vertices, where $B\subset\mathbb{R}^2$ is a fixed model triangle. For $j=1,\dots,r-1$, define $\sigma_j\subset\Delta$ to be
$\phi_\Delta^{-1}$ of the straight line segment in $B$ joining
$\phi_\Delta(p_j)$ to $\phi_\Delta(q_j)$. Being the pullback of a straight
(hence real-analytic) segment under a diffeomorphism, $\sigma_j$ is a
smooth simple arc from $p_j$ to $q_j$ with interior in
$\operatorname{int}(\Delta)$ (the segment's interior avoids $\partial B$
since $\phi_\Delta(p_j)$ and $\phi_\Delta(q_j)$ lie on different sides of
$B$ for $1\le j\le r-1$).

These arcs have pairwise disjoint interiors. Traversing $\partial\Delta$
once gives the cyclic order $p_0,p_1,\dots,p_{r-1},p_r(=q_r),q_{r-1},
\dots,q_1,q_0(=p_0)$, and $\phi_\Delta$ carries this to the same cyclic
order around $\partial B$. For $1\le j<k\le r-1$, the four points
$$\phi_\Delta(p_j),\phi_\Delta(p_k),\phi_\Delta(q_k),\phi_\Delta(q_j)$$ occur
in this cyclic order, so $\phi_\Delta(p_k)$ and $\phi_\Delta(q_k)$ both lie
on the same arc of $\partial B$ cut by $\{\phi_\Delta(p_j),\phi_\Delta(q_j)
\}$ — that is, the endpoint pairs are nested rather than interleaved. Since
$B$ is convex, two chords of a convex region with nested boundary endpoints
have disjoint interiors; pulling back through $\phi_\Delta^{-1}$, the same
holds for $\sigma_j$ and $\sigma_k$.

The arcs $\sigma_1,\dots,\sigma_{r-1}$, together with the two sides of
$\Delta$ other than $e_{\mathrm{out}}$ and $e_{\mathrm{out}}$ itself,
partition $\Delta$ into $r$ closed regions $P_1,\dots,P_r$, where $P_j$ is
bounded by $I_j$, $J_j$, and the arcs $\sigma_{j-1}$, $\sigma_j$ (with
$\sigma_0$, $\sigma_r$ understood as the vertices $p_0=q_0$ and $p_r=q_r$,
i.e.\ degenerate when $j=1$ or $j=r$). Each $P_j$ is a closed topological
disk, since it is bounded by a simple closed curve formed from finitely
many arcs; its boundary meets $\partial\Delta$ exactly in $I_j\cup J_j$,
since the only boundary segments used in $\partial P_j$ are
$I_j\subset\partial\Delta$ and $J_j\subset\partial\Delta$, the arcs
$\sigma_{j-1},\sigma_j$ lying in $\Delta$ with endpoints on $\partial\Delta$
but interior in $\operatorname{int}(\Delta)$. Thus each $P_j$ is a polygonal
disk, $\Delta=\bigcup_j P_j$, the $P_j$ have pairwise disjoint interiors,
and $P_j\cap\partial\Delta=I_j\cup J_j$.

By construction the order of the $J_j$ along $e_{\mathrm{out}}$ (via the
points $q_0,\dots,q_r$) agrees with the order of the $I_j$ along the other
two sides (via $p_0,\dots,p_r$), since both are indexed consistently by
$j=1,\dots,r$.
\end{proof}
\begin{proposition}[Sector decomposition]\label{prop:cut}
Let \(S\) be a smooth compact two-dimensional manifold, let
\(A\subset S\) be closed, and let \(G\) be a component of
\(S\setminus A\) with \(\partial G\neq\varnothing\). Then there exist a
countable locally finite smooth triangulation \(\mathcal T\) of \(G\), a
nonempty at most countable index set \(I\), pairwise disjoint nonempty
open sectors \(\{G_i\}_{i\in I}\), roots \(g_i\in\partial G\), and proper
polygonal rays \(\rho_i\subset G_i\) such that:
\begin{enumerate}
\item \(\operatorname{diam}(\Delta)
<\tfrac12\operatorname{dist}(\Delta,A)\) for every
\(\Delta\in\mathcal T\).
\item Each \(G_i\) is homeomorphic to the interior of a closed triangle.
\item \(g_i\in\operatorname{Acc}(\rho_i)\subset\partial G\) for every
\(i\in I\), and \(\{g_i:i\in I\}\) is dense in \(\partial G\).
\item For every \(i\in I\), there are compact polygonal disks
\(D_{i,1}\subset D_{i,2}\subset\cdots\subset G\) such that
\[
G_i=\bigcup_{m\ge1}\operatorname{int}_G(D_{i,m}),
\]
and \(G_i\setminus D_{i,m}\) is connected for every \(m\).
Set
\[
Q_{i,1}=G_i\cap D_{i,1},\qquad
Q_{i,m}=G_i\cap(D_{i,m}\setminus D_{i,m-1})
\quad(m\ge2).
\]
Then \(G_i=\bigsqcup_{m\ge1}Q_{i,m}\).
For \(m\ge2\), each \(Q_{i,m}\) lies in a triangle
\(\Delta_{i,m}\), and consecutive triangles are equal
or adjacent. The first layer \(Q_{i,1}\) lies in a finite
union of triangles which, together with \(\Delta_{i,2}\),
have a common vertex.
The ray \(\rho_i\) passes successively through the layers,
eventually leaving every \(D_{i,m}\).
\item
\(\bigcap_{m\geq1}\cl(\bigcup_{k\geq m}Q_{i,k})\subset\partial G\).
\item \(G\subset\bigcup_{i\in I}\cl(G_i)\), and
\(\{\cl(G_i)\}_{i\in I}\) is locally finite in \(G\).
\item Each \((\cl(G_i)\setminus G_i)\cap G\) is contained in a countable
locally finite union of piecewise smooth arcs and therefore has
\(\lambda\)-measure zero.
\end{enumerate}
\end{proposition}

\begin{proof}
By Lemma~\ref{lem:2d-accessible}, applied to the closed set
\(S\setminus G\), choose a nonempty at most countable dense set
\(Q=\{q_n:n\in N\}\subset\operatorname{AC}(G)\). Applying
Lemma~\ref{lem:whitehead-forest}, obtain a countable locally finite smooth
triangulation \(\mathcal T\) satisfying item~(1), pairwise vertex-disjoint
rays \(\{R_n\}_{n\in N}\subset D(\mathcal T)\) with
\(\operatorname{Acc}(|R_n|)=\{q_n\}\), and a spanning forest
\(F\subset D(\mathcal T)\) whose components
\(\{C_\alpha\}_{\alpha\in J}\) are special trees. Each prescribed ray
\(R_n\) is the distinguished ray of one component.

For \(\alpha\in J\), let \(H_\alpha\) be the union of the triangles whose
dual vertices lie in \(C_\alpha\). The sets \(H_\alpha\) cover \(G\).
Fix \(\alpha\), write the distinguished ray of \(C_\alpha\) as
\(V_0,V_1,\dots\), and decompose \(H_\alpha\) into sectors by propagating
colors toward this ray.

\noindent\emph{Coloring.}
Let \(e\) be a side of a triangle in \(H_\alpha\) not crossed by an edge of
\(C_\alpha\). If \(e\subset\partial H_\alpha\), assign it one new color.
If \(e\) is shared by two triangles of \(H_\alpha\), assign its two sides
— one per triangle — distinct new colors. By construction, colors are only
ever propagated across sides that \emph{are} crossed by an edge of
\(C_\alpha\).

Orient every edge of each finite tree \(F_\beta\subset C_\alpha\) toward
the distinguished ray, and process the triangles of \(H_\alpha\) in this
orientation, starting from the triangles corresponding to the leaves of
the \(F_\beta\) (the \emph{terminal triangles}, where all sides other than
the one crossed by the outgoing tree-edge are initially colored). For a
triangle \(\Delta\) being processed, let \(e_{\mathrm{out}}\) be the side
crossed by its outgoing edge; its other two sides are partitioned into
finitely many colored arcs, either assigned initially or propagated from
the triangles already processed.
Lemma~\ref{lem:triangle-propagation} then subdivides \(\Delta\) into
polygonal disks and propagates each incoming color onto a distinct
nondegenerate subarc of \(e_{\mathrm{out}}\).

Once all finite branches attached at \(V_k\) have been processed this way,
apply the same step to the triangle corresponding to \(V_k\) itself, with
the side crossed by \(V_kV_{k+1}\) as the outgoing side. At \(k=0\) there
is no incoming edge along the distinguished ray, but the remaining two
sides of the triangle at \(V_0\) are still partitioned into colored arcs
coming from its attached branches (or, if none, from the initial coloring
on its own non-\(C_\alpha\) sides).

Proceeding inductively along the distinguished ray colors every triangle
of \(H_\alpha\), each divided into finitely many polygonal pieces, as in
Figure~\ref{fig:process}.
\begin{figure}[h!]
\centering
\begin{tikzpicture}[scale=5.5]

    \coordinate (B) at (2, 0);
    \coordinate (C) at (1, 1.5);
    \coordinate (D) at (0.5, 0.75);
    \coordinate (E) at (1.5, 0.75);
    \coordinate (F) at (1, 0); 
    
    \coordinate (H) at (1.75, 0.375); 
    \coordinate (I) at (1, 0.75); 
    \coordinate (J) at (2.5, 1); 
  
    \coordinate (L) at (2.5, 1.5);
  \coordinate (Z) at ($(J)!0.5!(L)$); 
\coordinate (M_FH) at ($(F)!0.5!(H)$);
\coordinate (M_EJ) at ($(E)!0.5!(J)$);
 \coordinate (M_CI) at ($(C)!0.5!(I)$);
 \coordinate (G) at ($(D)!2/3!(M_CI)$);
\coordinate (M_CE) at ($(C)!0.5!(E)$);
 \coordinate (M) at ($(I)!2/3!(M_CE)$);
\coordinate (M_IF) at ($(I)!0.5!(F)$);
 \coordinate (N) at ($(D)!2/3!(M_IF)$);
\coordinate (M_EF) at ($(E)!0.5!(F)$);
 \coordinate (O) at ($(I)!2/3!(M_EF)$);
\coordinate (M_EH) at ($(E)!0.5!(H)$);
 \coordinate (P) at ($(F)!2/3!(M_EH)$);
\coordinate (M_HB) at ($(H)!0.5!(B)$);
 \coordinate (S) at ($(F)!2/3!(M_HB)$);
\coordinate (M_HJ) at ($(H)!0.5!(J)$);
 \coordinate (Q) at ($(E)!2/3!(M_HJ)$);
\coordinate (M_JL) at ($(J)!0.5!(L)$);
 \coordinate (R) at ($(E)!2/3!(M_JL)$);
 \fill[red, opacity=0.4] (D) -- (I) -- (M_CI) -- cycle;
 \coordinate (M_IE) at ($(I)!0.5!(E)$);
  \coordinate (M_CE) at ($(C)!0.5!(E)$);
 \coordinate (M_IM_IE) at ($(I)!0.5!(M_IE)$);
  \fill[red, opacity=0.4] (I) -- (M_IM_IE) -- (M_CI) -- cycle;
 \fill[blue, opacity=0.4] (D) -- (C) -- (M_CI) -- cycle;
  \fill[blue, opacity=0.4] (C) -- (M_CI) -- (M_IM_IE) -- (M_IE) -- cycle;
 \fill[green, opacity=0.4] (C) -- (M_IE) -- (E) -- cycle;
  \fill[yellow, opacity=0.4] (D) -- (M_IF) -- (I) -- cycle;
   \fill[blue, opacity=0.4] (D) -- (F) -- (M_IF) -- cycle;
 \coordinate (M_FM_EF) at ($(F)!0.5!(M_EF)$);
 \coordinate (M_EM_EF) at ($(E)!0.5!(M_EF)$); 
    \fill[blue, opacity=0.4] (F) -- (M_FM_EF) -- (M_IF) -- cycle;
\fill[yellow, opacity=0.4] (M_IF) -- (M_FM_EF) -- (M_EF) -- (I) -- cycle;
 \coordinate (M_M_EFE) at ($(M_EF)!0.5!(E)$); 
  \coordinate (M_M_M_EF) at ($(M_EF)!0.5!(M_M_EFE)$); 
  \fill[green, opacity=0.4] (E) -- (M_M_EFE) -- (M_IE) -- cycle;   
   \fill[blue, opacity=0.4] (M_M_M_EF) -- (M_M_EFE) -- (M_IE) -- (M_IM_IE) -- cycle; 
       \fill[red, opacity=0.4] (M_M_M_EF) -- (M_EF) -- (I) -- (M_IM_IE) -- cycle; 
  \fill[yellow, opacity=0.4] (B) -- (F) -- (M_FH) -- cycle;   
    \fill[blue, opacity=0.4] (B) -- (H) -- (M_FH) -- cycle; 
  \fill[yellow, opacity=0.4] (J) -- (H) -- (M_EJ) -- cycle;   
    \fill[blue, opacity=0.4] (L) -- (E) -- (Z) -- cycle;   
 \coordinate (M_M_EHE) at ($(M_EH)!0.5!(E)$);    
  \coordinate (M_EM_M_EHE) at ($(M_M_EHE)!0.5!(E)$);
 \coordinate (M_EM_EM_M_EHE) at ($(M_EM_M_EHE)!0.5!(E)$);
\fill[green, opacity=0.4] (E) -- (M_EM_M_EHE) -- (M_M_EFE) -- cycle; 
 \coordinate (M_M_EM_EM_M_EHE) at ($(M_EM_M_EHE)!0.5!(M_M_EHE)$);
 \fill[blue, opacity=0.4] (M_M_M_EF) -- (M_M_EFE) -- (M_EM_M_EHE) -- (M_M_EM_EM_M_EHE) -- cycle; 
 \fill[red, opacity=0.4] (M_M_M_EF) -- (M_EF)  -- (M_M_EHE) --(M_M_EM_EM_M_EHE) -- cycle;
  \coordinate (M_M_EM_M_EHE) at ($(M_M_EHE)!0.5!(M_EH)$);
 \fill[yellow, opacity=0.4] (M_M_EM_M_EHE) -- (M_FM_EF) -- (M_EF) -- (M_M_EHE) -- cycle;
  \fill[blue, opacity=0.4] (M_FM_EF) -- (F) -- (M_EH) -- (M_M_EM_M_EHE)--cycle;
  \coordinate (M_HM_EH) at ($(M_EH)!0.5!(H)$);    
  \fill[yellow, opacity=0.4] (M_FH) -- (F) -- (M_EH) -- (M_HM_EH) -- cycle;   
   \fill[blue, opacity=0.4] (M_HM_EH) -- (H) -- (M_FH) -- cycle;  

  \coordinate (M_M_EJ) at ($(M_EJ)!0.5!(E)$);
  \coordinate (M_M_M_EJ) at ($(E)!0.5!(M_M_EJ)$); 
\coordinate (M_M_M_M_EJ) at ($(E)!0.5!(M_M_M_EJ)$); 
\fill[green, opacity=0.4] (E) -- (M_EM_M_EHE) -- (M_M_M_M_EJ) -- cycle; 
\coordinate (M_M_M_M_M_EJ) at ($(M_M_M_EJ)!0.5!(M_M_M_M_EJ)$); 

 \fill[blue, opacity=0.4] (M_M_M_M_M_EJ) -- (M_M_M_M_EJ) -- (M_EM_M_EHE) -- (M_M_EM_EM_M_EHE) -- cycle; 
 \fill[red, opacity=0.4] (M_M_M_M_M_EJ) -- (M_M_M_EJ)  -- (M_M_EHE) --(M_M_EM_EM_M_EHE) -- cycle; 
 \coordinate (M_M_M_M_M_M_EJ) at ($(M_M_M_EJ)!0.5!(M_M_EJ)$); 
 \fill[yellow, opacity=0.4] (M_M_EM_M_EHE) -- (M_M_M_M_M_M_EJ) --  (M_M_M_EJ) -- (M_M_EHE) -- cycle;
   \fill[blue, opacity=0.4]  (M_M_M_M_M_M_EJ) -- (M_M_EJ)-- (M_EH) -- (M_M_EM_M_EHE)--cycle;
   \coordinate (MM_EJ) at ($(M_M_EJ)!0.5!(M_EJ)$); 
  \fill[yellow, opacity=0.4] (MM_EJ) --(M_M_EJ)-- (M_EH) -- (M_HM_EH) -- cycle;  
    \fill[blue, opacity=0.4] (M_HM_EH) -- (H) -- (M_EJ) --(MM_EJ) -- cycle; 
  \coordinate (M_ZJ) at ($(Z)!0.5!(J)$);
    \coordinate (M_M_ZJ) at ($(M_ZJ)!0.5!(Z)$);
      \coordinate (M_M_M_ZJ) at ($(M_M_ZJ)!0.5!(M_ZJ)$);
 \fill[yellow, opacity=0.4] (J) -- (M_ZJ) -- (M_EJ) -- cycle;  
     \fill[blue, opacity=0.4] (M_M_M_ZJ) -- (M_ZJ)-- (M_EJ) --(MM_EJ) -- cycle; 
 \fill[yellow, opacity=0.4] (MM_EJ) --(M_M_EJ)-- (M_M_ZJ) -- (M_M_M_ZJ) -- cycle;  
 \coordinate (ZZ) at ($(Z)!0.5!(M_M_ZJ)$);
  \coordinate (ZZZ) at ($(ZZ)!0.5!(M_M_ZJ)$);
    \fill[blue, opacity=0.4]  (M_M_M_M_M_M_EJ) -- (M_M_EJ)-- (M_M_ZJ) -- (ZZZ)--cycle; 
\fill[yellow, opacity=0.4] (ZZZ) -- (M_M_M_M_M_M_EJ) --  (M_M_M_EJ) -- (ZZ) -- cycle;  
\coordinate (ZZZZ) at ($(Z)!0.5!(ZZ)$);
\fill[green, opacity=0.4] (E) -- (M_M_M_M_EJ) -- (ZZZZ) --(Z)-- cycle; 
\coordinate (X) at ($(ZZZZ)!0.5!(ZZ)$);
 \fill[blue, opacity=0.4] (M_M_M_M_M_EJ) -- (M_M_M_M_EJ) -- (ZZZZ) -- (X) -- cycle; 
 \fill[red, opacity=0.4] (M_M_M_M_M_EJ) -- (M_M_M_EJ)  -- (ZZ) --(X) -- cycle; 
   \draw[thick] (D) -- (C);
    \draw[thick] (I) -- (E);
   \draw[thick] (C) -- (E);
   \draw[thick] (F) -- (B);
   \draw[thick] (H) -- (J);
   \draw[thick] (D) -- (I);
    \draw[thick] (F) -- (C); 
    \draw[thick] (F) -- (D); 
    \draw[thick] (F) -- (E); 

    \draw[thick] (F) -- (H); 
    \draw[thick] (E) -- (H) -- (B); 
    \draw[thick] (E) -- (J); 
      \draw[thick] (E) -- (L) -- (J);
\fill[black] (G) circle (0.02); 
\fill[black] (M) circle (0.02); 
\fill[black] (N) circle (0.02); 
\fill[black] (O) circle (0.02); 
\fill[black] (P) circle (0.02); 
\fill[black] (S) circle (0.02); 
\fill[black] (Q) circle (0.02); 
\fill[black] (R) circle (0.02); 
\draw[very thick] (G) -- (M_CI);
\draw[very thick] (M) -- (M_CI);
\draw[very thick] (M) -- (M_IE);
\draw[very thick] (O) -- (M_IE);
\draw[very thick] (N) -- (M_IF);
\draw[very thick] (O) -- (M_IF);
\draw[very thick] (O) -- (M_EF);
\draw[very thick] (P) -- (M_EF);
\draw[very thick] (P) -- (M_EH);
\draw[very thick] (S) -- (M_FH);
\draw[very thick] (P) -- (M_FH);
\draw[very thick] (M_EH) -- (Q);
\draw[very thick] (Q) -- (M_EJ);
\draw[very thick] (R) -- (M_EJ);
\draw[very thick] (Z) -- (R);
\draw[->, very thick] (Z) -- ++(0.3, -0.3) node[right] {$C_{\alpha}$};

\end{tikzpicture}
\caption{Coloring \(H_\alpha\). The special tree \(C_\alpha\) is drawn in bold; its vertices are the barycenters of the triangles in \(H_\alpha\).}
\label{fig:process}
\end{figure}
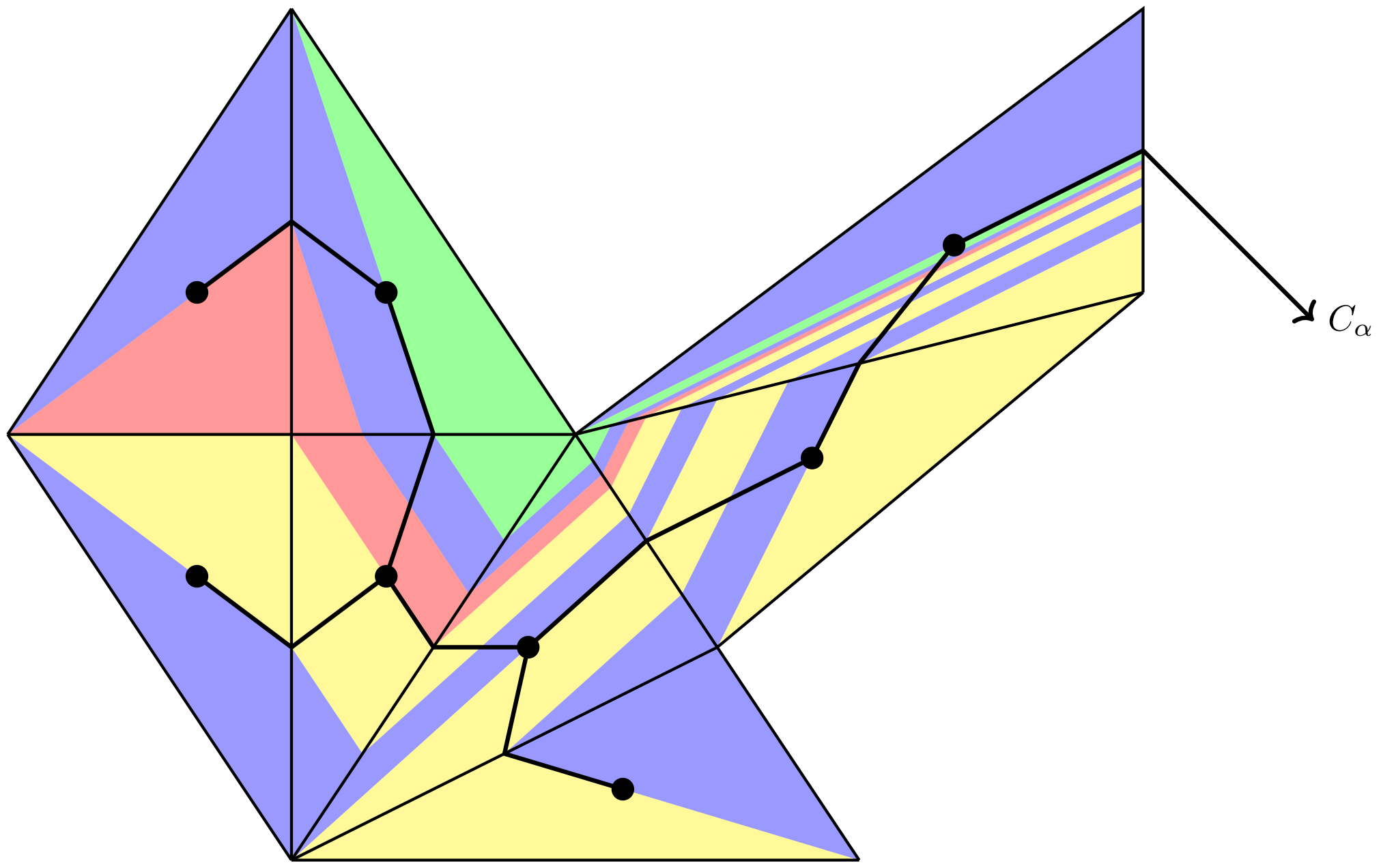

\medskip
Fix a color \(\xi\). By construction, every piece of color \(\xi\) has
exactly one successor and, except for the initial piece, exactly one
predecessor. Enumerate these pieces as
\(P^\alpha_{\xi,1},P^\alpha_{\xi,2},\ldots\), where consecutive pieces
meet along a nondegenerate boundary arc.
Nonconsecutive pieces may meet only at triangulation vertices. By
Lemma~\ref{lem:triangle-propagation}, every triangulation vertex of a propagated
piece belongs to its incoming arc; hence every such intersection
involving \(P^\alpha_{\xi,l}\) is contained in
\(P^\alpha_{\xi,l}\cap P^\alpha_{\xi,l-1}\).

Consequently, for every \(k<l\),
\[
P^\alpha_{\xi,l}\cap
\bigl(P^\alpha_{\xi,k}\cup\cdots\cup P^\alpha_{\xi,l-1}\bigr)
=
P^\alpha_{\xi,l}\cap P^\alpha_{\xi,l-1},
\]
and this intersection is a nondegenerate boundary arc of \(P^\alpha_{\xi,l}\).
Every point of the arc is a limit of points of
\(\operatorname{int}_G(P^\alpha_{\xi,l})\), which is disjoint from
\(P^\alpha_{\xi,k}\cup\cdots\cup P^\alpha_{\xi,l-1}\).
Hence the arc also lies in the boundary of that union. Hence
Lemma~\ref{lem:gluing}, applied inductively, shows that every finite
union \(P^\alpha_{\xi,k}\cup\cdots\cup P^\alpha_{\xi,l}\) is a compact
polygonal disk. Set
\[
D^\alpha_{\xi,m}=\bigcup_{j=1}^m P^\alpha_{\xi,j},
\qquad
G^\alpha_\xi=
\operatorname{int}_G\Bigl(\bigcup_{m\ge1}D^\alpha_{\xi,m}\Bigr).
\]

\smallskip\noindent\emph{Local finiteness (proves item~(6)).}
Fix \(x\in G\) and choose a relatively compact neighbourhood \(W\) of
\(x\). Since \(\mathcal T\) is locally finite, \(\cl(W)\) meets only
finitely many triangles, each containing finitely many colored pieces. If
\(\cl(G_\xi^\alpha)\) meets \(W\), then \(G_\xi^\alpha\) meets \(W\), so
one of its colored pieces lies in one of these finitely many triangles.
Hence only finitely many sector closures meet \(W\), and
\(\{\cl(G_\xi^\alpha)\}_{\alpha,\xi}\) is locally finite in \(G\).

\smallskip\noindent\emph{Each sector is a disk (proves item~(2)).}
By definition,
\(G_\xi^\alpha=\operatorname{int}_G(\bigcup_{m\geq1}D_{\xi,m}^\alpha)\).
If \(x\in G_\xi^\alpha\), choose a neighbourhood \(W\) of \(x\) contained
in \(\bigcup_{m\geq1}D_{\xi,m}^\alpha\). By local finiteness, \(W\) meets
only finitely many colored pieces, and hence
\(W\subset D_{\xi,m}^\alpha\) for some \(m\). Thus
\(G_\xi^\alpha=\bigcup_{m\geq1}\operatorname{int}_G(D_{\xi,m}^\alpha)\).
Consequently, every compact subset of \(G_\xi^\alpha\) lies in some
\(\operatorname{int}_G(D_{\xi,m}^\alpha)\), so every loop is
null-homotopic. Hence \(G_\xi^\alpha\) is a connected simply connected
noncompact surface, and therefore
\(G_\xi^\alpha\cong\mathbb R^2\cong\operatorname{int}(\mathbb B)\).

\smallskip\noindent
\emph{Layers and their properties (proves item~(4)).}
Fix \(\alpha,\xi\), and suppress these indices temporarily.
Write \(C_m=P_m\cap P_{m+1}\) for the attaching arcs.
The arc \(C_1\) contains exactly one triangulation vertex
\(v\). By the propagation rule, every triangulation vertex
of \(C_m\), for \(m\ge2\), also belongs to \(C_{m-1}\).
Thus \(v\) is the only triangulation vertex occurring in
these arcs, and the arcs containing it form an initial
segment \(C_1,\ldots,C_k\). This segment is finite because
the color follows a simple path in the dual tree and only
finitely many triangles contain \(v\).

Replace \(P_1,\ldots,P_k\) by their union and relabel
the pieces:
\[
\widetilde P_1=P_1\cup\cdots\cup P_k,\qquad
\widetilde P_m=P_{k+m-1}\quad(m\ge2).
\]
The new successive attaching arcs are disjoint: any
intersection between the incoming and outgoing arcs of
a piece must be a triangulation vertex, and \(v\) no
longer occurs after the first attaching arc.
All original triangles containing the merged pieces,
as well as the triangle containing \(\widetilde P_2\),
contain \(v\).

We use again the notation \(P_{\xi,m}^\alpha\) for the regrouped pieces and set
\(
D_{\xi,m}^\alpha=\bigcup_{j=1}^{m}P_{\xi,j}^\alpha.
\)
The regrouping does not change the union of the disks, hence we still have
\(
G_\xi^\alpha
=\bigcup_{m\ge1}\operatorname{int}_G(D_{\xi,m}^\alpha).
\)

Define
\[
Q_{\xi,1}^\alpha=G_\xi^\alpha\cap D_{\xi,1}^\alpha,
\qquad
Q_{\xi,m}^\alpha
=G_\xi^\alpha\cap
(D_{\xi,m}^\alpha\setminus D_{\xi,m-1}^\alpha)
\quad(m\ge2).
\]
Then \(G_\xi^\alpha=\bigsqcup_{m\ge1}Q_{\xi,m}^\alpha\).
For \(m\ge2\), the layer \(Q_{\xi,m}^\alpha\) lies in
the triangle \(\Delta_{\xi,m}^\alpha\) containing
\(P_{\xi,m}^\alpha\), and consecutive triangles share
an edge. The first layer lies in a finite union of
triangles which, together with \(\Delta_{\xi,2}^\alpha\),
have the common vertex \(v\).

For every \(m\), the set
\(G_\xi^\alpha\setminus D_{\xi,m}^\alpha\) is connected:
the remaining pieces form a chain, and their interiors
join across the relative interiors of the attaching arcs.

Choose a point \(p_m\) in the relative interior of each
arc \(P_{\xi,m}^\alpha\cap P_{\xi,m+1}^\alpha\).
Inside \(P_{\xi,m+1}^\alpha\), join \(p_m\) to \(p_{m+1}\)
by a simple polygonal arc whose interior lies in
\(\operatorname{int}_G(P_{\xi,m+1}^\alpha)\).
These arcs concatenate to an embedded polygonal ray
\(\rho_\xi^\alpha\subset G_\xi^\alpha\), starting at
\(p_1\) and passing successively through the layers.
Every compact subset of \(G\) meets only finitely many
colored pieces, so the ray eventually leaves every
such compact set. Hence \(\rho_\xi^\alpha\) is proper
in \(G\).

\smallskip\noindent\emph{The root (proves item~(3), accessibility part).}
If the distinguished ray of \(C_\alpha\) is the prescribed ray \(R_j\),
then the pieces \(P_{\xi,m}^\alpha\) eventually follow the same sequence
of shrinking triangles as \(R_j\). Hence
\(\operatorname{Acc}(\rho_\xi^\alpha)=\{q_j\}\), and we set
\(g_\xi^\alpha=q_j\).

For every other sector, compactness of \(S\) and properness of
\(\rho_\xi^\alpha\) in \(G\) imply that
\(\operatorname{Acc}(\rho_\xi^\alpha)\) is a nonempty subset of
\(\partial G\). Choose any \(g_\xi^\alpha\) in this set. Thus, in all
cases,
\(g_\xi^\alpha\in\operatorname{Acc}(\rho_\xi^\alpha)\subset\partial G\).

\smallskip\noindent\emph{High layers accumulate only on \(\partial G\)
(proves item~(5)).}
Let \(y_l\in Q_{\xi,k_l}^\alpha\) with \(k_l\to\infty\). Since the
corresponding triangles \(\Delta_{\xi,k_l}^\alpha\) leave every compact
subset of \(G\), the sequence \((y_l)\) has no accumulation point in
\(G\). As \(S\) is compact, all its accumulation points therefore lie in
\(\partial G\). Hence
\[
\bigcap_{m\geq1}\cl\Bigl(\bigcup_{k\geq m}Q_{\xi,k}^\alpha\Bigr)
\subset\partial G.
\]

\smallskip\noindent\emph{Measure-zero relative boundary (proves item~(7)).}
The interior gluing arcs disappear from the boundary of the union. Hence
\((\cl(G_\xi^\alpha)\setminus G_\xi^\alpha)\cap G\) is contained in the
initial boundary arc together with the two outer boundary chains of the
pieces \(P_{\xi,m}^\alpha\). By
Lemma~\ref{lem:triangle-propagation}, these chains are contained in a
locally finite union of polygonal arcs, which has \(\lambda\)-measure
zero. Thus
\(
\lambda\bigl((\cl(G_\xi^\alpha)\setminus G_\xi^\alpha)\cap G\bigr)=0.
\)

\smallskip\noindent\emph{Covering (completes item~(6)).}
Every colored piece \(P_{\xi,m}^\alpha\) has interior contained in
\(G_\xi^\alpha\), and therefore
\(P_{\xi,m}^\alpha\subset\cl(G_\xi^\alpha)\). Since the colored pieces
cover every triangle of \(H_\alpha\), we have
\(H_\alpha\subset\bigcup_\xi\cl(G_\xi^\alpha)\). As the sets
\(H_\alpha\) cover \(G\), the sector closures cover \(G\).

\smallskip\noindent\emph{Relabeling and density of the roots (completes
item~(3)).}
Relabel the nonempty at most countable family
\(\{G_\xi^\alpha\}_{\alpha,\xi}\) as \(\{G_i\}_{i\in I}\), and relabel
the associated roots, rays, disks, and layers accordingly. If \(R_j\) is
the distinguished ray of \(C_\alpha\), then every sector following
\(R_j\) has root \(q_j\); hence each \(q_j\) occurs among the roots
\(g_i\). Since \(Q=\{q_j:j\in N\}\) is dense in \(\partial G\), so is
\(\{g_i:i\in I\}\). This completes item~(3).

\end{proof}

Let \(\mathbb B\subset\mathbb R^2\) be the closed triangle
with vertices \(w_1=(0,0)\), \(w_2=(1,0)\), and
\(w_3=(\tfrac12,1)\). For \(s,t\in[0,1]\), put
\[
\gamma_s(t)=(1-t)(s,0)+tw_3.
\]
Thus \(\gamma_s\) parametrizes the segment from \((s,0)\)
to \(w_3\).

Choose \(0=t_0<t_1<t_2<\cdots\) with \(t_m\to1\), and
partition \(\operatorname{int}(\mathbb B)\) into the layers
\[
L_m=\{\gamma_s(t):0<s<1,\ t_{m-1}<t\le t_m\},
\qquad m\ge1.
\]
Write
\(
\ell_{w_3}=\{\gamma_{1/2}(t):0<t<1\}
\)
for the central fiber, and
\(
\ell^+=\{\gamma_{1/2}(t):t_1\le t<1\}
\)
for its terminal ray.

\begin{lemma}[Straightening a sector]\label{lem:straightening}
Let \(G_i\) be a sector with disks \(D_{i,m}\), layers
\(Q_{i,m}\), and ray \(\rho_i\) constructed in
Proposition~\ref{prop:cut}. Then there exists a homeomorphism
\(
\phi_i:\operatorname{int}(\mathbb B)\to G_i
\)
such that
\[
\phi_i(L_m)=Q_{i,m}\quad(m\ge1),
\qquad
\phi_i(\ell^+)=\rho_i.
\]
Moreover, for every \(m\ge1\), the restriction of \(\phi_i\)
to \(L_m\cup L_{m+1}\) extends continuously to its closure
in \(\mathbb B\), with values in \(\cl_S(G_i)\).
\end{lemma}

\begin{proof}
Fix \(i\) and suppress this index. Let
\[
C_m=D_m\cap P_{m+1}=P_m\cap P_{m+1},
\]
and let \(p_m\) be the point where \(\rho\) crosses \(C_m\).
By construction, \(C_m\) and \(C_{m+1}\) are disjoint
boundary arcs of \(P_{m+1}\), and the subarc of \(\rho\)
from \(p_m\) to \(p_{m+1}\) has interior in
\(\operatorname{int}_G(P_{m+1})\).

Set
\[
E_m=\{\gamma_s(t):0\le s\le1,\ 0\le t\le t_m\}.
\]
We construct compatible homeomorphisms
\(
\phi_m:E_m\to D_m
\)
such that the terminal arc of \(E_m\) is mapped onto \(C_m\),
\(
\phi_m(\gamma_{1/2}(t_m))=p_m,
\)
and the central segment
\[
\{\gamma_{1/2}(t):t_1\le t\le t_m\}
\]
is mapped onto the subarc of \(\rho\) from \(p_1\) to \(p_m\).

For \(m=1\), choose a boundary homeomorphism
\(\partial E_1\to\partial D_1\) taking the terminal arc
onto \(C_1\) and \(\gamma_{1/2}(t_1)\) to \(p_1\), and
extend it over the disks.

Suppose that \(\phi_m\) has been constructed. The disk
\(E_{m+1}\) is obtained from \(E_m\) by attaching
\[
\{\gamma_s(t):0\le s\le1,\ t_m\le t\le t_{m+1}\},
\]
while \(D_{m+1}\) is obtained from \(D_m\) by attaching
\(P_{m+1}\) along \(C_m\). By Lemma~\ref{lem:gluing},
\(\phi_m\) extends to a homeomorphism
\(
\phi_{m+1}:E_{m+1}\to D_{m+1}
\)
which maps the new terminal arc onto \(C_{m+1}\) and the
new central segment onto the subarc of \(\rho\) from
\(p_m\) to \(p_{m+1}\).

The compatible maps define
\(
\phi_i:\operatorname{int}(\mathbb B)\to G_i.
\)
Since
\[
\operatorname{int}(\mathbb B)
=\bigcup_{m\ge1}\operatorname{int}(E_m),
\qquad
G_i=\bigcup_{m\ge1}\operatorname{int}_G(D_m),
\]
every point of either space lies in the interior of some
finite-stage disk. There the map and its inverse agree with
finite-stage homeomorphisms, so \(\phi_i\) is a
homeomorphism.

Fix \(m\) and let \(z\in\operatorname{int}(B)\). Choose \(k\ge m\)
with \(z\in\operatorname{int}(E_k)\). Since \(\varphi_k\) is injective
and \(\varphi_k(E_m)=\varphi_m(E_m)=D_m\), we have
\(\varphi_i(z)=\varphi_k(z)\in D_m\) if and only if \(z\in E_m\).
Thus \(\varphi_i(E_m\cap\operatorname{int}(B))=D_m\cap G_i\),
and the definitions of the layers give \(\varphi_i(L_m)=Q_{i,m}\).
The matching of the central segments gives
\(
\phi_i(\ell^+)=\rho_i.
\)
Finally,
\(\overline{L_m\cup L_{m+1}}\subset E_{m+1}\), so
\(\phi_{m+1}\) gives the required continuous extension.
\end{proof}

\begin{lemma}\label{lem:function}
There exists a homeomorphism $f_{\mathbb B}:\operatorname{int}(\mathbb B)
\to\operatorname{int}(\mathbb B)$ such that:
\begin{enumerate}
\item For every $x\in\operatorname{int}(\mathbb B)$,
$f_{\mathbb B}^{\,n}(x)\to w_3$.
\item For every $x\in\operatorname{int}(\mathbb B)$,
$\operatorname{dist}(f_{\mathbb B}^{\,n}(x),\ell_{w_3})\to0$; moreover, the
forward orbit of $x$ meets every sufficiently high layer.
\item For every $m\ge1$, $f_{\mathbb B}(L_m)\subset L_m\cup L_{m+1}$.
\item For every $m\ge1$, $f_{\mathbb B}^{-1}(L_m)\subset L_{m-1}\cup L_m$,
where $L_0=\varnothing$.
\item If $y_n\in\operatorname{int}(\mathbb B)$ and $y_n\to\partial\mathbb
B$, then $d_{\mathbb B}(f_{\mathbb B}(y_n),y_n)\to0$ and $d_{\mathbb
B}(f_{\mathbb B}^{-1}(y_n),y_n)\to0$.
\end{enumerate}
\end{lemma}

\begin{proof}
Every point of $\operatorname{int}(\mathbb B)$ can be written uniquely as
$x(s,t)=(1-t)(s,0)+tw_3$ for $0<s<1$, $0<t<1$. Choose a continuous function
$\lambda:(0,1)\to(0,1/4]$ with $\lambda(s)\to0$ as $s\to0$ or $s\to1$. For
$m\geq1$, set $v_m=(1-\lambda(s))t_m+\lambda(s)t_{m+1}$, and put $v_0=0$.

For each fixed $s$, let $\rho_s:(0,1)\to(0,1)$ be affine on every interval
$[t_m,t_{m+1}]$ with $\rho_s(t_m)=v_m$. Since $t_m<v_m<t_{m+1}$ and
$v_m\to1$, $\rho_s$ is an increasing homeomorphism of $(0,1)$ onto itself,
and $\rho_s(t)>t$ for every $0<t<1$. Define $f_{\mathbb B}(x(s,t))=
x(s,\rho_s(t))$. Since $\rho_s(t)$ depends continuously on $(s,t)$ and the
affine pieces agree at layer boundaries, $f_{\mathbb B}$ is a homeomorphism,
with inverse $f_{\mathbb B}^{-1}(x(s,u))=x(s,\rho_s^{-1}(u))$.

\smallskip\noindent\emph{Convergence to the apex (item~1).} Fix
$x=x(s,t)$ and write $f_{\mathbb B}^{\,n}(x)=x(s,t_n)$. Since
$t_{n+1}=\rho_s(t_n)>t_n$, the sequence $(t_n)$ increases to some
$\ell\leq1$; if $\ell<1$, continuity would give $\rho_s(\ell)=\ell$, a
contradiction. Thus $t_n\to1$, so $f_{\mathbb B}^{\,n}(x)\to w_3$.

\smallskip\noindent\emph{Convergence to the central fiber and layer
visits (item~2).} Since $\operatorname{dist}(x(s,t_n),\ell_{w_3})\leq
\|x(s,t_n)-w_3\|\to0$, the first claim of item~2 follows. Since
$\rho_s([t_{m-1},t_m])=[v_{m-1},v_m]\subset[t_{m-1},t_{m+1}]$, we get
$f_{\mathbb B}(L_m)\subset L_m\cup L_{m+1}$ (item~3) and similarly
$f_{\mathbb B}^{-1}(L_m)\subset L_{m-1}\cup L_m$ (item~4). Since $t_n\to1$
and one iterate cannot skip a layer, every forward orbit meets every
sufficiently high layer, completing item~2.

\smallskip\noindent\emph{Boundary displacement (item~5).} Since $x(s,u)-
x(s,t)=(u-t)(w_3-(s,0))$, there is $C>0$ with $d_{\mathbb B}(x(s,u),x(s,t))
\leq C|u-t|$. For $t\in[t_m,t_{m+1}]$, $|\rho_s(t)-t|\leq\lambda(s)\max\{
t_{m+1}-t_m,t_{m+2}-t_{m+1}\}$, which tends to $0$ uniformly as the layer
index $m\to\infty$, since $t_{m+1}-t_m\to0$.

On the first layer, $\rho_s(t)=\tfrac{v_1}{t_1}t$, so $|\rho_s(t)-t|=
\lambda(s)\tfrac{t_2-t_1}{t_1}t$; choosing $\lambda$ small enough that
$\lambda(s)\tfrac{t_2-t_1}{t_1}\leq1$ gives $|\rho_s(t)-t|\leq t$, so the
displacement tends to $0$ as $t\to0$. Finally $|\rho_s(t)-t|\leq\lambda(s)$,
so the displacement tends to $0$ uniformly as $s\to0$ or $s\to1$. Combining
these three regimes (large layer index, $t\to0$, and $s\to0,1$) — which
together cover every way a point can approach $\partial\mathbb B$ — shows
that if $y_n\to\partial\mathbb B$ then $d_{\mathbb B}(f_{\mathbb B}(y_n),
y_n)\to0$.

For the inverse, suppose \(y_n\to\partial B\) and set
\(z_n=f_B^{-1}(y_n)\). Then \(z_n\to\partial B\): otherwise some
compact set \(K\subset\operatorname{int}(B)\) would contain infinitely
many \(z_n\), and the compact set \(f_B(K)\subset\operatorname{int}(B)\)
would contain the corresponding \(y_n\). Applying the forward
displacement estimate to \(z_n\) gives
\[
   d_B(f_B^{-1}(y_n),y_n)=d_B(z_n,f_B(z_n))\longrightarrow 0.
\]
Thus $f_{\mathbb B}$ has all the required properties.
\end{proof}

\begin{proposition}\label{prop:verify-axioms}
Let \(S\) be a smooth compact two-dimensional manifold, let
\(A\subset S\) be nonempty and closed, and let \(G\) be a component of
\(S\setminus A\) with \(\partial G\neq\varnothing\). Let
\(\{G_i,g_i,\rho_i\}_{i\in I}\) be the sectors, roots, and rays given by
Proposition~\ref{prop:cut}. For each \(i\in I\), let
\(\phi_i:\operatorname{int}(\mathbb B)\to G_i\)
be the homeomorphism given by
Lemma~\ref{lem:straightening}, and define
\(f_i:\cl(G_i)\to\cl(G_i)\) by
\[
f_i(x)=
\begin{cases}
\phi_i\circ f_{\mathbb B}\circ\phi_i^{-1}(x),&x\in G_i,\\
x,&x\in\cl(G_i)\setminus G_i.
\end{cases}
\]
Then \(\{(G_i,g_i,f_i)\}_{i\in I}\) is a homeomorphic
root--contracting decomposition of \(G\) toward \(A\). Moreover,
\[
d(f_i^{\pm1}(x),x)<\tfrac54\,d(x,A)
\]
for every \(i\in I\) and \(x\in G_i\).
\end{proposition}

\begin{proof}
We verify the conditions of Definition~\ref{def:sector-decomposition}.

\smallskip\noindent\emph{Sector structure.}
By Proposition~\ref{prop:cut}, the sets $G_i$ are nonempty, open, and
pairwise disjoint, and $G\subset\bigcup_i\cl(G_i)$. Moreover, $g_i\in
\cl(G_i)\cap\partial G\subset A$, and the roots are dense in $\partial G$.
By definition, $f_i$ fixes $\cl(G_i)\setminus G_i$ pointwise; for every $i$,
let $\phi_i:\operatorname{int}(\mathbb B)\to G_i$ be the homeomorphism given
by Lemma~\ref{lem:straightening}.

\smallskip\noindent\emph{The maps on the sector closures.}
On $G_i$, $f_i$ is conjugate to $f_{\mathbb B}$; it remains to prove
continuity at $\cl(G_i)\setminus G_i$. Let $x_n\in G_i$, $x_n\to x\in
\cl(G_i)\setminus G_i$, and put $y_n=\phi_i^{-1}(x_n)$. Then $y_n$ leaves
every compact subset of $\operatorname{int}(\mathbb B)$, hence approaches
$\partial\mathbb B$, so by Lemma~\ref{lem:function},
$d_{\mathbb B}(f_{\mathbb B}(y_n),y_n)\to0$.

Let $m_n$ be the layer index containing $y_n$; passing to a subsequence,
assume $(m_n)$ bounded or $m_n\to\infty$. If $(m_n)$ is bounded, we may
assume all $y_n$ lie in one fixed layer $L_m$; since $f_{\mathbb
B}(L_m)\subset L_m\cup L_{m+1}$, both $y_n$ and $f_{\mathbb B}(y_n)$ lie in
$L_m\cup L_{m+1}$, and by choice of $\phi_i$ its restriction there extends
continuously to the closure, giving $d_S(f_i(x_n),x_n)\to0$. If $m_n\to
\infty$, then $x_n$ and $f_i(x_n)$ lie in the same polygonal piece or in
two consecutive pieces, contained in one triangle of $\mathcal T$ or two
sharing an edge; since their layer indices $\to\infty$, these triangles
leave every compact subset of $G$ with diameters $\to0$, again giving
$d_S(f_i(x_n),x_n)\to0$. Thus $f_i(x_n)\to x=f_i(x)$, so $f_i$ is
continuous on $\cl(G_i)$.

Define $f_i^{-1}(x)=\phi_i\circ f_{\mathbb B}^{-1}\circ\phi_i^{-1}(x)$ on
$G_i$ and $=x$ on $\cl(G_i)\setminus G_i$; the same argument, using
$f_{\mathbb B}^{-1}(L_m)\subset L_{m-1}\cup L_m$, shows $f_i^{-1}$
continuous on $\cl(G_i)$. Hence $f_i$ is a homeomorphism of $\cl(G_i)$ onto
itself.

\smallskip\noindent\emph{The rooted $\omega$-limit condition.}
Fix $x\in G_i$, $y=\phi_i^{-1}(x)$. By Lemma~\ref{lem:function}, $f_{\mathbb
B}^{\,n}(y)\to w_3$, so the orbit of $x$ eventually leaves every compact
subset of $G_i$. If $p\in\omega_{f_i}(x)$, choose $n_j\to\infty$ with
$f_i^{\,n_j}(x)\to p$; the layer indices of $f_{\mathbb B}^{\,n_j}(y)$ tend
to infinity, and by Proposition~\ref{prop:cut}(5) high layers accumulate
only on $\partial G$, so $p\in\partial G\subset A$. Thus
$\omega_{f_i}(x)\subset A$.

For $g_i\in\omega_{f_i}(x)$: since
$g_i\in\operatorname{Acc}(\rho_i)$ and
$\phi_i(\ell^+)=\rho_i$, there are points $z_m\in\ell^+$
in layers $L_{k_m}$ with $k_m\ge2$, $k_m\to\infty$,
and $\phi_i(z_m)\to g_i$. By
Lemma~\ref{lem:function}, the orbit of $y$ meets every sufficiently high
layer, so we may choose $n_m\to\infty$ with $f_{\mathbb B}^{\,n_m}(y)\in
L_{k_m}$, hence $f_i^{\,n_m}(x),\phi_i(z_m)\in Q_{i,k_m}\subset
\Delta_{i,k_m}\in\mathcal T$ with $\operatorname{diam}(\Delta_{i,k_m})\to0$.
Thus $d_S(f_i^{\,n_m}(x),\phi_i(z_m))\to0$, and since $\phi_i(z_m)\to g_i$,
$f_i^{\,n_m}(x)\to g_i$, so $g_i\in\omega_{f_i}(x)$.

\smallskip\noindent\emph{Slowing down near $A$.}
Fix $b\in A$ and $\varepsilon>0$. If $b\notin\partial G$, then
$d_G:=\operatorname{dist}(b,G)>0$, and
$V=B(b,d_G/2)$ is disjoint from $G$, so the claim is vacuous.

Assume now that \(b\in\partial G\), and let \(x\in G_i\).
Choose original triangles \(\Delta,\Delta'\in\mathcal T\)
containing \(x\) and \(f_i(x)\), respectively, with
\(\Delta\cap\Delta'\neq\varnothing\).
Such a choice is possible because away from the first
layer the triangles are equal or adjacent, while the
triangles forming the first layer and the triangle
containing the second layer have a common vertex. Since
\[
\operatorname{diam}(\Delta)
 <\tfrac12\operatorname{dist}(\Delta,A)
 \le\tfrac12 d_S(x,A),
\]
the first case gives
$d_S(f_i(x),x)\le\tfrac12d_S(x,A)$. In the second case, choose
$z\in\Delta\cap\Delta'$. Then
\[
\operatorname{dist}(\Delta',A)
 \le d_S(z,A)
 \le d_S(x,A)+\operatorname{diam}(\Delta)
 <\tfrac32d_S(x,A),
\]
and hence
$$
d_S(f_i(x),x)
 \le\operatorname{diam}(\Delta)+\operatorname{diam}(\Delta')
 <\tfrac12d_S(x,A)+\tfrac34d_S(x,A)
 =\tfrac54d_S(x,A).
$$
The same argument applies to \(f_i^{-1}\), since it also
moves points only between the same or neighboring layers,
and the corresponding triangles can again be chosen
to intersect. Thus
$d_S(f_i^{\pm1}(x),x)<\tfrac54d_S(x,A)$ for every $x\in G_i$.

Set $V=B(b,4\varepsilon/5)$. If $x\in G_i\cap V$, then
$d_S(x,A)\le d_S(x,b)<4\varepsilon/5$, and therefore
$d_S(f_i^{\pm1}(x),x)<\varepsilon$.

\smallskip\noindent\emph{Local finiteness.}
The family $\{\cl(G_i)\}_{i\in I}$ is locally finite in $G$ by
Proposition~\ref{prop:cut}(6).

Thus \(\{(G_i,g_i,f_i)\}_{i\in I}\) satisfies the componentwise
root--contracting conditions on \(G\). Moreover, the displacement estimate
is independent of the choice of the component \(G\).
\end{proof}

\begin{proof}[Proof of Theorem~\ref{thm:2D}]
If \(A=S\), the identity map realizes \(A\). Hence assume
\(A\subsetneq S\). Assume first that $S$ is connected, and let $A\subset S$ be nonempty and
closed. For each component \(G_\alpha\) of \(S\setminus A\), apply
Proposition~\ref{prop:verify-axioms} and denote the resulting data by
\(
\{(G_{\alpha,i},g_{\alpha,i},f_{\alpha,i})\}_{i\in I_\alpha}.
\)
The uniform displacement estimate shows that the slowing-down condition
holds simultaneously over all components. Local finiteness away from
\(A\) follows because every point of \(S\setminus A\) has a neighbourhood
contained in its unique complementary component. Hence
\(
\{(G_{\alpha,i},g_{\alpha,i},f_{\alpha,i})
  :\alpha\in\Theta,\ i\in I_\alpha\}
\)
is a homeomorphic root--contracting decomposition of
\(S\setminus A\) toward \(A\).

The union of the roots is dense in \(\partial A\). Indeed, let
\(a\in\partial A\) and let \(U\) be a neighbourhood of \(a\). Choose
\(r>0\) such that \(B(a,r)\subset U\), and take
\(x\in B(a,r/3)\setminus A\). Let \(G_\alpha\) be the component of
\(S\setminus A\) containing \(x\). Since \(a\notin G_\alpha\),
\[
\operatorname{dist}(x,\partial G_\alpha)
\le d(x,a)<r/3.
\]
Hence \(\partial G_\alpha\cap B(a,r)\neq\varnothing\). Since the roots are
dense in \(\partial G_\alpha\), the neighbourhood \(U\) contains a root.

Proposition~\ref{prop:cut}(7) gives
\(
\lambda\bigl((\cl(G_{\alpha,i})\setminus G_{\alpha,i})\setminus A\bigr)=0
\)
for every sector. Therefore
Theorem~\ref{thm:root-decomp-implies-attractor} yields a homeomorphism
$f:S\to S$ with $\lambda(B(A))=1$, and
Corollary~\ref{cor:gsp-implies-maf} gives $\Lambda(f)=A$.

Now suppose that $S$ has connected components $M_1,\dots,M_r$. If
$A\cap M_j\neq\varnothing$ for every $j$, apply the connected case to
each pair $(M_j,A\cap M_j)$ and combine the resulting homeomorphisms.

If \(A\) misses some component, choose \(p\in A\). On each component
meeting \(A\), use the map supplied by the connected case, and map each
component disjoint from \(A\) constantly to \(p\). Since the components of
\(S\) are clopen, these maps combine into a continuous self-map
\(f:S\to S\). Every point in a component disjoint from \(A\) enters \(A\)
after one iteration, while \(A\) has full basin in every component meeting
it; hence \(\lambda(B(A))=1\).

To prove minimality, let \(B\subsetneq A\) be closed. Then
\(B\cap M_j\subsetneq A\cap M_j\) for some component \(M_j\) meeting
\(A\). By the connected case, \(A\cap M_j\) is the Milnor global
attractor of \(f|_{M_j}\), so \(B\cap M_j\) does not have full basin in
\(M_j\). Hence \(B\) does not have full basin in \(S\).

Now let \(K\subset S\) be any closed set with full basin. Then
\(A\cap K\) also has full basin, since almost every \(\omega\)-limit set
is contained in both \(A\) and \(K\). By the preceding argument,
\(A\cap K\) cannot be a proper subset of \(A\). Thus \(A\subset K\), so
\(A\) is the smallest closed set with full basin and therefore
\(\Lambda(f)=A\).

Thus every smooth compact two-dimensional manifold is an MAF--space,
and the realizing map may be chosen to be a homeomorphism whenever
\(A\) meets every connected component of \(S\).
\end{proof}

\begin{cor}
Let \(S\) be a closed connected topological surface equipped with a
nonatomic Borel probability measure \(\mu\) that is positive on every
nonempty open set. Then every nonempty closed subset of \(S\) is the
Milnor global attractor of a homeomorphism of \(S\).
\end{cor}

\begin{proof}
Every topological surface admits a smooth structure; see
Hatcher~\cite[Theorem A]{HatcherSurfaces}. Choose such a structure and let \(\lambda\) be a
smooth area probability measure. By the Oxtoby--Ulam theorem \cite{OxtobyUlam}, there is
a homeomorphism \(h:S\to S\) such that
\(
h_*\mu=\lambda.
\)

Let \(A\subset S\) be nonempty and closed. By
Theorem~\ref{thm:2D}, there is a homeomorphism \(F:S\to S\) whose
Milnor global attractor with respect to \(\lambda\) is \(h(A)\). Set
\[
f=h^{-1}\circ F\circ h.
\]
For every closed \(C\subset S\),
\[
h(B_f(C))=B_F(h(C)).
\]
Since \(h_*\mu=\lambda\),
\[
\mu(B_f(C))=\lambda(B_F(h(C))).
\]
Thus \(C\) has full basin for \(f\) if and only if \(h(C)\) has
full basin for \(F\). In particular, \(A\) has full basin. Moreover,
if \(C\subset S\) is closed and has full basin for \(f\), then
\(h(C)\) has full basin for \(F\). Since \(h(A)\) is the smallest
closed full-basin set for \(F\),
\(
h(A)\subset h(C),
\)
and hence \(A\subset C\). Therefore \(A\) is the smallest closed
full-basin set for \(f\), and so it is the Milnor global attractor
of \(f\).
\end{proof}

\begin{cor}
Let \(S\) be a smooth compact connected two-dimensional manifold and
let \(A\subset S\) be nonempty and closed. The homeomorphism in
Theorem~\ref{thm:2D} may be chosen so that
\[
\Lambda(f)=\Gamma(f)=A.
\]
\end{cor}

\begin{proof}
Let \(f\) be the homeomorphism constructed in the proof of
Theorem~\ref{thm:2D}. By Proposition~\ref{prop:cut}(7), the
union of the sector boundaries outside \(A\) is meagre. Every other
point lies either in \(A\) or in the interior of a sector, and hence
belongs to \(B_f(A)\). Thus \(B_f(A)\) is residual.

Let \(C\subsetneq A\) be closed and choose \(a\in A\setminus C\) and
an open neighbourhood \(U\subset S\setminus C\) of \(a\). If
\(U\cap\operatorname{int}(A)\neq\varnothing\), then \(U\) contains a
nonempty open set on which \(\omega_f(x)=\{x\}\not\subset C\).
Otherwise \(a\in\partial A\), and density of the roots gives a root
\(g_{\alpha,i}\in U\); then \(U\cap G_{\alpha,i}\) is nonempty and open,
and \(g_{\alpha,i}\in\omega_f(x)\not\subset C\) for every
\(x\in U\cap G_{\alpha,i}\). Hence \(B_f(C)\) is not residual.
Therefore \(\Gamma(f)=A\), while \(\Lambda(f)=A\) by
Theorem~\ref{thm:2D}.
\end{proof}

\section{Beyond manifolds}\label{sec:beyond}

The root--contracting mechanism of Section~\ref{sec:root-contracting}
relies on two features of a manifold: enough room in each complementary
component to build sectors, and enough self-maps to realize the required
local contractions. Away from manifolds, both can fail, and the MAF
property can fail with them. We collect a few examples illustrating this,
together with a list of open questions.

\subsection{A rigid space: the Cook continuum}

\begin{exmp}[Cook continuum \cite{Cook}]\label{ex:cook}
There exists a nondegenerate compact connected metric space $K$ -- the
\emph{Cook continuum} -- every continuous self-map of which is either the
identity or a constant map.
\end{exmp}

\begin{proposition}
The Cook continuum $K$ is not an MAF--space, for any Borel probability
measure on $K$ with full support.
\end{proposition}

\begin{proof}
If $f\equiv x_0$, then $\Lambda(f)=\{x_0\}$; if $f=\operatorname{id}_K$,
then $\Lambda(f)=K$. Hence $\{\Lambda(f):f:K\to K\text{ continuous}\}
=\{\{x_0\}:x_0\in K\}\cup\{K\}$, which omits every proper nondegenerate
subcontinuum of $K$. Thus $K$ is not an MAF--space.
\end{proof}

This argument uses no special feature of $K$: any \emph{rigid continuum}
(every continuous self-map is the identity or constant) fails to be an
MAF--space.

\subsection{A branching obstruction: the tripod}

\begin{exmp}[Tripod]\label{ex:tripod}
Let $G$ be a tripod: three arcs $[v,a_1],[v,a_2],[v,a_3]$ joined only at
$v$, with free endpoints $a_1,a_2,a_3$.
\end{exmp}

\begin{proposition}\label{prop:tripod}
Let \(A=\{a_1\}\). Then \(A\) is not the Milnor global attractor of any
homeomorphism of \(G\), with respect to \(H^1\).
\end{proposition}

\begin{proof}
Any homeomorphism \(f:G\to G\) fixes \(v\) and permutes the three open
legs \(U_i=\operatorname{int}([v,a_i])\). Hence there exists \(p\geq1\)
such that \(f^p(U_3)=U_3\). For every \(x\in U_3\), we therefore have
\(f^{kp}(x)\in U_3\) for all \(k\geq0\).

Since \(\cl(U_3)=[v,a_3]\) is compact and disjoint from \(A\), the
sequence \((f^{kp}(x))_{k\geq0}\) has an accumulation point in
\(\cl(U_3)\subset G\setminus A\). Thus
\(\omega_f(x)\not\subset A\), so \(U_3\subset G\setminus B(A)\).
Because \(H^1(U_3)>0\), the basin \(B(A)\) does not have full measure.
Therefore \(A\) cannot be the Milnor global attractor of \(f\).
\end{proof}

The constant map \(x\mapsto a_1\) realizes the singleton \(\{a_1\}\)
as its Milnor global attractor. Thus some prescribed closed sets on
the tripod require a noninvertible realizing map. For a general
branched continuum \(X\), it remains open to characterize
\(\{\Lambda(h):h\in\operatorname{Homeo}(X)\}\).

\subsection{A measure-theoretic obstruction: the prolonged sine curve}

\begin{exmp}[Prolonged sine curve]\label{ex:sine-curve}
Let $\Sigma=\{(x,\sin(1/x)):0<x\le1\}$ be a sine tail, and let
$I=\{0\}\times[-1,1]$, so that $\cl(\Sigma)\setminus\Sigma=I$. Let
$J=\{0\}\times[-1,2]$, the segment $I$ extended past its endpoint
$(0,1)$ up to $(0,2)$, and set $X=J\cup\Sigma$; thus $J$ and $\Sigma$ are
the two path components of $X$, with $J\setminus I=\{0\}\times(1,2]$.
\end{exmp}

\begin{proposition}\label{prop:sine-not-maf}
Let $F=\{a,b\}$ with $a\in J\setminus I$, $b\in\Sigma$, and let $\mu$ be
any reference measure with $\mu(I)>0$. Then $F$ is not the Milnor global
attractor of any continuous $g:X\to X$, but $F$ \emph{is} realized as a
generic limit set by some continuous $f:X\to X$.
\end{proposition}

\begin{proof}
\emph{Generic realization.}
Parametrize the sine tail by
\(\gamma(t)=(1/t,\sin t)\), \(t\in[t_0,\infty)\), where \(t_0=1\), choose \(t_b\) with
\(\gamma(t_b)=b\), and let \(\varphi:[t_0,\infty)\to[t_0,\infty)\) be a
continuous map fixing \(t_b\), satisfying
\(\varphi^n(t)\to t_b\) for every \(t\), and
\(|\varphi(t)-t|\to0\) as \(t\to\infty\). Defining
\(f(\gamma(t))=\gamma(\varphi(t))\) makes every point of \(\Sigma\)
converge to \(b\), while the vanishing-displacement condition implies
that \(f\) extends continuously over \(I\) by the identity. On \(J\),
choose a continuous interval map fixing \(I\) pointwise and attracting
every point of \(J\setminus I\) to \(a\). Thus
\(\omega_f(x)=\{a\}\) on \(J\setminus I\) and
\(\omega_f(x)=\{b\}\) on \(\Sigma\). Since \(I\) is nowhere dense,
\(X\setminus I\) is residual, so \(B(F)\) is residual. Neither
\(B(\{a\})\) nor \(B(\{b\})\) is residual, since they miss the nonempty
open sets \(\Sigma\) and \(J\setminus I\), respectively. Hence
\(\Gamma(f)=F\).

\smallskip
\noindent\emph{Not Milnor-realizable.}
Suppose \(\Lambda(g)=F=\{a,b\}\). For almost every \(x\),
\(\omega_g(x)\subset F\). Since
\(g(\omega_g(x))=\omega_g(x)\), we have
\(\omega_g(x)\subset g(F)\). Thus \(g(F)\) has full basin, and
minimality gives \(F\subset g(F)\). Since \(F\) has two points,
\(g(F)=F\), so \(g|_F\) is a permutation. Suppose
first that \(g\) swaps \(a\) and \(b\). The sets \(J\) and \(\Sigma\)
are the two path components of \(X\). Since continuous images of
path-connected sets are path-connected, \(g(a)=b\in\Sigma\) implies
\(g(J)\subset\Sigma\), while \(g(b)=a\in J\) implies
\(g(\Sigma)\subset J\).

Now let \(y\in I\) and choose \(x_n\in\Sigma\) with \(x_n\to y\).
Since \(g(x_n)\in J\) and \(J\) is closed, continuity gives
\(g(y)\in J\). On the other hand, \(y\in J\) and
\(g(J)\subset\Sigma\), so \(g(y)\in\Sigma\), contradicting
\(J\cap\Sigma=\varnothing\). Hence \(a\) and \(b\) cannot be swapped
and must therefore be fixed. Thus
\(a\) and \(b\) are fixed, so
\(g(J)\subset J\) and \(g(\Sigma)\subset\Sigma\). Repeating the same
limit argument shows \(g(I)\subset I\). Since \(I\cap F=\varnothing\),
every orbit starting in \(I\) remains in \(I\), so
\(I\subset X\setminus B(F)\). As \(\mu(I)>0\), this contradicts
\(\Lambda(g)=F\).
\end{proof}

Thus Milnor and generic realization are genuinely different notions on the
same space, and this example also gives a simple instance of a continuum
that is not an MAF--space: since $F$ cannot be realized as a Milnor global
attractor of any continuous self-map of $X$, the space $X$ itself fails
the MAF property.

\subsection{A non-locally-connected example: root--contracting is not
necessary}

\begin{exmp}[Spiral compactum]\label{ex:spiral}
Let \(C=S^1\) and define
\[
\gamma(t)=\bigl((1-t)\cos(1/t),(1-t)\sin(1/t)\bigr),
\qquad t\in(0,1],
\]
and set \(\Sigma=\gamma((0,1])\) and \(X=C\cup\Sigma\). Thus \(\Sigma\)
is a spiral accumulating on \(C\).
\end{exmp}
The space \(X\) is not locally connected at the points of \(C\), so it
falls outside the hypotheses of
Theorem~\ref{thm:root-decomp-implies-attractor}. Nevertheless, it shows
that, for general compact metric spaces, the existence of a
root--contracting decomposition is not necessary for realizability.

Suppose, toward a contradiction, that \(\Sigma=X\setminus C\) admits a
root--contracting decomposition toward \(C\). We will show that this
forces the root of the unique sector to be dense in \(C\), contradicting
the definition of a root--contracting decomposition.

Let \(U\subset\Sigma\) be a sector with associated map \(h\). Since the
root of \(U\) belongs to \(\cl(U)\cap C\), this intersection is
nonempty.

We first show that some component of \(U\) has closure meeting \(C\). Suppose
otherwise. Then \(U\) has infinitely many components accumulating on
\(C\). At most one contains the terminal point \(\gamma(1)\), so choose
a component \(J\) that does not. Its closure is a compact interval in
\(\Sigma\), with endpoints \(p,q\in\cl(U)\setminus U\), both fixed by
\(h\).

Since \(h(\cl(J))\) is path-connected and contains \(p\) and \(q\), it
lies in \(\Sigma\). Identifying \(\Sigma\) with \((0,1]\), write
\(\cl(J)=[t_p,t_q]\). As \(h(p)=p\) and \(h(q)=q\),
\(h(\cl(J))\supset[t_p,t_q]\). If
\(h(\cl(J))\not\subset\cl(J)\), then some \(\zeta\in J\) satisfies
\(h(\zeta)\notin[t_p,t_q]\). Assume \(h(\zeta)>t_q\), the other case
being symmetric. By the Intermediate Value Theorem, there exists
\(x\in J\) with \(h(x)=q\). Since \(q\) is fixed,
\(\omega_h(x)\cap C=\varnothing\), a contradiction. Hence
\(h(\cl(J))\subset\cl(J)\), so
\(\omega_h(x)\subset\cl(J)\subset X\setminus C\) for every \(x\in J\),
again a contradiction.

Thus some component of \(U\) has closure meeting \(C\), and every such
component contains a tail \(\gamma((0,c))\). Consequently every sector
contains a tail. Since any two tails intersect, there is at most one
sector, whose unique root therefore cannot be dense in \(C\).

On the other hand, \(C\) is directly realizable. Choose
\(\theta>0\) with \(\theta/(2\pi)\notin\mathbb Q\), let
\(R_\theta:C\to C\) be rotation by \(\theta\), define
\(\tau(t)=t/(1+\theta t)\), and set
\(f|_C=R_\theta\) and
\(f(\gamma(t))=\gamma(\tau(t))\).
Since \(1/\tau(t)=1/t+\theta\), this extends continuously to \(R_\theta\)
on \(C\). Moreover, \(\tau^n(t)=t/(1+n\theta t)\to0\), while the angular
coordinate advances by \(n\theta\). As \(\theta/(2\pi)\) is irrational,
\(\omega_f(x)=C\) for every \(x\in X\). Thus \(C\) is the Milnor global
attractor of \(f\), although \(X\setminus C\) admits no
root--contracting decomposition toward \(C\).

\subsection{Open questions}

\begin{itemize}

\item Characterize the one-dimensional continua that are MAF--spaces.
How does the answer change when one restricts from continuous maps to
homeomorphisms?
\item For a compact metric space \(X\), characterize
\[
\mathcal H(X)=
\{\Lambda(h):h\in\operatorname{Homeo}(X)\}.
\]
Equivalently, which closed sets can occur as Milnor global attractors of
homeomorphisms of \(X\)?

\item Given a closed set \(A\subset X\) and a continuous map
\(g:A\to A\), when does there exist a continuous map \(f:X\to X\)
extending \(g\) such that \(A\) is the Milnor global attractor of \(f\)?

\item For which closed sets \(A\subset X\) does there exist a continuous
map \(f:X\to X\) such that \(A\) is the Milnor global attractor of \(f\)
and
\(
\omega_f(x)=A
\)
for almost every \(x\in X\)?

\item What is the relationship between Milnor and generic realization?
Under what hypotheses does
\(
\Lambda(f)=\Gamma(f)?
\)

\item Does the MAF property extend to smooth compact manifolds of dimension
three or higher?

\item To what extent does realizability depend on the reference measure?
What changes for measures that are not mutually absolutely continuous
with a smooth volume measure?

\item Can the realizing maps on smooth manifolds be chosen to be
diffeomorphisms, or more generally \(C^r\) maps?

\end{itemize}

\section*{Acknowledgements}
The authors were supported by institutional support for the development
of research organizations in the Czech Republic (RVO) for IČ~47813059.

The authors used AI to help identify gaps in the proofs and check the validity of individual steps. The mathematical ideas are the authors’ own, and the authors take full responsibility for the correctness of the proofs.

\end{document}